\documentclass[11pt]{amsart}
\usepackage[a4paper,margin=27mm]{geometry}
\usepackage{amssymb,mathtools,booktabs}
\usepackage[hidelinks]{hyperref}
\hypersetup{
  pdftitle={Integer group determinants for the metacyclic group C(p squared) semidirect C(p)},
  pdfauthor={Chatchawan Panraksa}
}
\theoremstyle{plain}
\newtheorem{theorem}{Theorem}[section]
\newtheorem{lemma}[theorem]{Lemma}
\newtheorem{proposition}[theorem]{Proposition}
\newtheorem{corollary}[theorem]{Corollary}
\newtheorem{mainthm}{Theorem}

\theoremstyle{definition}
\newtheorem{problem}[theorem]{Problem}
\theoremstyle{remark}
\newtheorem{remark}[theorem]{Remark}
\newtheorem{example}[theorem]{Example}
\numberwithin{equation}{section}

\newcommand{\Z}{\mathbf Z}
\newcommand{\Q}{\mathbf Q}
\newcommand{\F}{\mathbf F}
\newcommand{\cO}{\mathcal O}
\newcommand{\cN}{\mathcal N}
\newcommand{\cB}{\mathcal B}
\newcommand{\cL}{\mathcal L}
\newcommand{\p}{\mathfrak p}
\newcommand{\PP}{\mathfrak P}
\newcommand{\m}{\mathfrak m}
\newcommand{\Lam}{\Lambda}
\newcommand{\Fbar}{\overline F}
\newcommand{\Ftil}{\widetilde F}
\newcommand{\Gp}{G_p}
\newcommand{\trz}{\operatorname{tr}_0}
\DeclareMathOperator{\nrd}{nrd}
\DeclareMathOperator{\Tr}{Tr}
\DeclareMathOperator{\Nm}{N}
\DeclareMathOperator{\tr}{tr}
\DeclareMathOperator{\adj}{adj}
\DeclareMathOperator{\SL}{SL}
\DeclareMathOperator{\M}{M}
\DeclareMathOperator{\diag}{diag}
\DeclareMathOperator{\Gal}{Gal}
\DeclareMathOperator{\ord}{ord}

\begin{document}

\title[Integer group determinants for $C_{p^2}\rtimes C_p$]{Integer group determinants for the metacyclic group $C_{p^2}\rtimes C_p$}

\author{Chatchawan Panraksa}
\thanks{Corresponding author: Chatchawan Panraksa
(\texttt{chatchawan.pan@mahidol.ac.th}). ORCID:
\href{https://orcid.org/0000-0003-0692-5453}{0000-0003-0692-5453}.}
\address{Applied Mathematics Program, Mahidol University International College,
Salaya, Nakhon Pathom 73170, Thailand}
\email{chatchawan.pan@mahidol.ac.th}


\begin{abstract}
Let $p$ be an odd prime and let $G_p=C_{p^2}\rtimes C_p$ be the nonabelian group
of order $p^3$ and exponent $p^2$. We study the integers that occur as
determinants of integral group matrices of $G_p$. For every odd $p$ we show that a
determinant prime to $p$ is congruent modulo $p^3$ to the $p$th power of a circulant
determinant coming from the quotient of $G_p$ by its elements of order dividing $p$;
that every integer $m$ with $m^{p-1}\equiv1\pmod{p^3}$ occurs; and that a
determinant divisible by $p$ is divisible by $p^{2p+3}$, which is sharp. For $p=3$
and $p=5$ we determine the set of determinants completely. For $p=3$, combining
this classification with the known cases gives a strict inclusion chain for
the determinant sets of all five groups of order $27$. For $p=5$ the answer involves
a cubic invariant of the prime ideals of the fifth cyclotomic field; one consequence is that $5^{13}$ is not a
determinant, although determinants of $5$-adic valuation $13$ exist.
\end{abstract}

\maketitle

\section{Introduction}\label{sec:intro}

Let $G$ be a finite group. To an element $F=\sum_{g\in G}a_gg$ of the integral group
ring $\Z[G]$ we attach the integral group matrix $(a_{gh^{-1}})_{g,h\in G}$. Its
determinant
\[
 D_G(F)=\det\bigl(a_{gh^{-1}}\bigr)_{g,h\in G}
\]
is the determinant of left multiplication by $F$ on $\Z[G]$, so $D_G$ is
multiplicative. Olga Taussky-Todd~\cite{TausskyTodd1977} asked which integers arise
in this way, that is, for a description of the set
\[
 S(G)=\{D_G(F):F\in\Z[G]\}
\]
of \emph{integer group determinants} of $G$. The cyclic case was taken up by
Laquer~\cite{Laquer1980} and Newman~\cite{Newman1980b,Newman1980}. Mahoney and
Newman~\cite{MahoneyNewman1980} showed that an abelian group of order $n$ has every
integer prime to $n$ as a group determinant exactly when it is cyclic, and
Mahoney~\cite{Mahoney1982} studied some nonabelian groups with two generators.
Pinner and Smyth~\cite{PinnerSmyth2020} determined $S(G)$ for every group of order
at most $14$. There
is also a link with Lehmer's problem. For abelian $G$, Lind's logarithmic Mahler
measure is $m_G(F)=|G|^{-1}\log|D_G(F)|$~\cite{Lind2005}, and the same formula is used
for nonabelian groups~\cite{MossinghoffPinnerHeisenberg}. The least value of
$|D_G(F)|$ exceeding $1$ therefore gives the Lind--Lehmer constant
\[
 \lambda(G)=|G|^{-1}\log\min\{|D_G(F)|:F\in\Z[G],\ |D_G(F)|\ge2\}.
\]

Now let $p$ be an odd prime. Of the five groups of order $p^3$, four have been
studied in detail. Following Mossinghoff and Pinner, we write $S^*(G)$ for the
elements of $S(G)$ prime to $p$ and $S_0(G)$ for those divisible by $p$, and we put
\[
 T_k=\{m\in\Z: m^{p-1}\equiv1\pmod{p^k}\}\qquad(k\ge1).
\]
Every integer prime to $p$ lies in $S(C_{p^3})$ (Laquer~\cite{Laquer1980} and
Newman~\cite{Newman1980b}). For the elementary abelian group, $S^*(C_p^3)=T_3$ by work
of DeSilva and Pinner~\cite[Lemmas~2.1 and~2.2]{DeSilvaPinner2014}; Trafford~\cite[Ch.~2,
Theorem~12]{Trafford1992} attributes the criterion for all ranks to
Odoni~\cite{Odoni1983}. For the Heisenberg group $H_p$, $S^*(H_p)=T_3$
\cite[Corollary~3.2]{MossinghoffPinnerHeisenberg}. Mossinghoff and Pinner also
proved that $S_0(H_p)\subseteq p^{p^2+3}\Z$, with the exponent sharp, and that
$S(H_3)=\{m:m\equiv\pm1\pmod{27}\}\cup3^{12}\Z$~\cite[Theorems~3.4
and~3.6]{MossinghoffPinnerHeisenberg}. More recently they showed that
$S_0(C_p\times C_{p^2})\subseteq p^{2p+3}\Z$, again with the exponent sharp, and
determined $S^*(C_3\times C_9)$~\cite[Theorems~3.2 and~3.4]{MossinghoffPinner2026}.
Their subsequent classification of $S_0(C_3\times C_9)$ completes the abelian
cases of order $27$~\cite[Theorem~1]{MossinghoffPinnerAbelian27}.

The subject of this paper is the fifth group,
\[
 \Gp=\langle x,y\mid x^{p^2}=y^p=1,\ yxy^{-1}=x^{1+p}\rangle\cong C_{p^2}\rtimes C_p ,
\]
the nonabelian group of order $p^3$ and exponent $p^2$. Mossinghoff and Pinner
identified this family as a remaining case in~\cite[\S1]{MossinghoffPinnerHeisenberg}
and singled out $G_3$ at the end of~\cite[\S1]{MossinghoffPinnerAbelian27}.
Mahoney's divisibility theorems for $Z$-metacyclic
groups~\cite[Theorems~7--9]{Mahoney1982} do not apply to $\Gp$: in his terminology a
$Z$-metacyclic group has cyclic abelianization, whereas $\Gp/\Gp'\cong C_p\times
C_p$. Nor is $\Gp$ among the semidirect products $C_q\rtimes C_n$, $q$ prime,
studied by Paudel and Pinner~\cite{PaudelPinner2026Semidirect}.

\subsection*{Results for every odd prime}

The elements of order dividing $p$ in $\Gp$ form a normal subgroup
$N=\langle x^p,y\rangle\cong C_p\times C_p$, and $\Gp/N$ is cyclic of order $p$,
generated by the image $X$ of $x$. For $F\in\Z[\Gp]$ let $\Ftil\in\Z[C_p]$ be the
image of $F$ under $x\mapsto X$, $y\mapsto1$.

\begin{mainthm}\label{thm:A}
Let $p$ be an odd prime and $F\in\Z[\Gp]$.
\begin{enumerate}
\item[(a)] If $p\nmid D_{\Gp}(F)$, then
\[
 D_{\Gp}(F)\equiv D_{C_p}(\Ftil)^p\pmod{p^3}.
\]
\item[(b)] Every integer $m$ with $m^{p-1}\equiv1\pmod{p^3}$ lies in $S(\Gp)$; that
is, $T_3\subseteq S^*(\Gp)$.
\item[(c)] If $p\mid D_{\Gp}(F)$, then $p^{2p+3}\mid D_{\Gp}(F)$. The exponent is
attained: $p^{2p+3}$ exactly divides $D_{\Gp}(p+(1-x)+(1-y)^2)$ for $p\ge5$, and
$D_{G_3}(x^4+x^7+x^6y^2)=3^9$.
\end{enumerate}
\end{mainthm}

The $p$th powers in $(\Z/p^3\Z)^\times$ are exactly the classes of $T_2$, so
part~(a) implies $S^*(\Gp)\subseteq T_2$. That inclusion also follows from
$S^*(N)=T_2$ and the containment theorem of Yamaguchi and
Yamaguchi~\cite[Theorem~1.4]{YamaguchiDedekind2024}. But part~(a) says more: modulo
$p^3$, the determinant depends only on the image of $F$ in $\Gp/N\cong C_p$, whose
own group determinants take every value prime to $p$. For $C_p\times C_{p^2}$ the same congruence can be deduced from a theorem of
De Silva, Mossinghoff, Pigno and Pinner~\cite{DSMPP2019} (Remark~\ref{rem:siblings});
for the Heisenberg group, in which every element has order dividing $p$, it fails.

Part~(a), combined with the Frobenius factorization (Lemma~\ref{lem:fact}), gives a
further necessary condition. To state it, let $\omega$ be a primitive $p$th root of
unity; this notation is kept throughout. A positive
integer $n$ is a \emph{norm from $\Z[\omega]$} if $n=\Nm_{\Q(\omega)/\Q}(\beta)$ for
some $\beta\in\Z[\omega]$, and a \emph{norm divisor} of an integer $m$ if moreover
$n\mid m$. A norm from $\Z[\omega]$ is the norm of an ideal, so a prime $q\ne p$
divides it to an exponent divisible by the order of $q$ modulo $p$; when $\Q(\omega)$
has class number one, as it does for $p\le19$~\cite[Theorem~11.1]{Washington1997},
the converse holds. In Proposition~\ref{prop:nec} we show that every
$m\in S^*(\Gp)$ lies in $T_2$ and has a norm divisor $n$ with $mn^p\in T_3$. A
condition of the same shape appears in Mossinghoff and Pinner's description of
$S^*(C_3\times C_9)$~\cite[Theorem~3.4]{MossinghoffPinner2026}. Taking $n=1$ gives
exactly the values in part~(b).

Parts~(b) and~(c) and the norm-divisor condition together decide which primes occur.

\begin{corollary}\label{cor:primes}
A prime $q$ belongs to $S(\Gp)$ if and only if $q^{p-1}\equiv1\pmod{p^3}$. Moreover
$2^p\in S(\Gp)$, so $\lambda(\Gp)\le p^{-2}\log2$, with equality for $p=3$ and
$p=5$: $\lambda(G_3)=\frac{1}{27}\log8$ and $\lambda(G_5)=\frac1{125}\log32$.
\end{corollary}

For an abelian $p$-group $C_p\times H$ of order $p^n$, Yamaguchi and
Yamaguchi~\cite[Theorem~1]{YamaguchiYamaguchi2023} proved that a prime $q$ is a group
determinant if and only if $q^{p-1}\equiv1\pmod{p^n}$; their Corollary~4, which holds
for all $p$-groups, gives only $q\in T_2$ for $\Gp$. For $p=7$ the inequality is
strict: $18^3\equiv1\pmod{343}$, so $18\in T_3\subseteq S(G_7)$, and $18<2^7$.

\subsection*{The case \texorpdfstring{$p=3$}{p=3}}

For $p=3$ the necessary condition is also sufficient, and part~(c) of
Theorem~\ref{thm:A} is the only restriction on the multiples of $3$. Since
$\Z[\omega]$ is now a principal ideal domain, a positive integer prime to $3$ is a
norm from $\Z[\omega]$ exactly when every prime $q\equiv2\pmod3$ divides it to an
even power.

\begin{mainthm}\label{thm:B}
Let $G_3=C_9\rtimes C_3$.
\begin{enumerate}
\item[(a)] An integer $m$ with $3\nmid m$ lies in $S(G_3)$ if and only if
$m\equiv\pm1\pmod 9$ and $m$ has a norm divisor $n$ with $mn^3\equiv\pm1\pmod{27}$.
\item[(b)] $S_0(G_3)=3^9\Z$. Consequently the $3$-adic valuation of a nonzero
element of $S(G_3)$ is either $0$ or at least $9$, and every such valuation occurs.
\end{enumerate}
\end{mainthm}

The divisor $n=1$ in~(a) gives $T_3=\{m:m\equiv\pm1\pmod{27}\}$, and larger divisors
give further values such as $62=2\cdot31$ (with $n=31$) and $91=7\cdot13$ (with
$n=7$). Neither of these is a group determinant of $C_3\times C_9$.
To prove part~(b), we give five explicit elements whose nontrivial representation
factors have the same product $3^8$. Adding multiples of the sum of all group
elements then gives every multiple of $3^9$. The two parts together give
$S(C_3\times C_9)\subsetneq S(G_3)$ and place $S(G_3)$ in a strict inclusion chain
with the determinant sets of the other four groups of order $27$
(Corollary~\ref{cor:p3}).

\subsection*{The case \texorpdfstring{$p=5$}{p=5}}

For $p=5$ the norm-divisor condition is not sufficient: $418=2\cdot11\cdot19$
satisfies it with $n=11$, but it is not a group determinant of $G_5$
(Proposition~\ref{prop:418}). The obstruction is a cubic invariant of the primes
above $11$, which we now define. Let $\cO=\Z[\omega]$, with $\omega$ a primitive
fifth root of unity. Since $\Phi_5(z)\equiv(z-1)^4\pmod5$, we may identify
$\cO/5\cO$ with $\F_5[w]/(w^4)$ by sending $\omega$ to the truncated exponential
\[
 E(w)=1+w+\frac{w^2}{2}+\frac{w^3}{6}.
\]
For $\alpha\in\cO$ prime to $5$ with image $\alpha(w)$, let $\ell_3(\alpha)\in\F_5$
be the coefficient of $w^3$ in $\log\bigl(\alpha(w)/\alpha(0)\bigr)$, where
$\log(1+z)=z-z^2/2+z^3/3$. Let $\cN$ be the set of positive integers of the form
$\Nm(\alpha)$ with $\alpha\in\cO$ prime to $5$ and $\ell_3(\alpha)=0$, and let $\cB$
be the set of primes $q\equiv1\pmod5$ such that $\ell_3(\gamma)\ne0$ for a generator
$\gamma$ of a prime ideal above $q$; by Section~\ref{sec:g5-prime-test}, the choice
of the prime ideal and of its generator does not matter. The smallest elements of
$\cB$ are $11$, $31$, $41$, $61$ and $71$.

\begin{mainthm}\label{thm:C}
Let $G_5=C_{25}\rtimes C_5$.
\begin{enumerate}
\item[(a)] An integer $m$ with $5\nmid m$ lies in $S(G_5)$ if and only if
$m^4\equiv1\pmod{25}$ and some $n\in\cN$ dividing $|m|$ satisfies
$(mn^5)^4\equiv1\pmod{125}$.
\item[(b)] We have
\[
 S_0(G_5)=5^{14}\Z\ \cup\ \bigcup_{q\in\cB}5^{13}q\Z .
\]
In particular, $5^{13}\notin S(G_5)$, the least positive element of $S_0(G_5)$ is
$5^{14}$, and the least positive group determinant of $5$-adic valuation $13$ is
$11\cdot5^{13}$.
\end{enumerate}
\end{mainthm}

Since every norm prime to $5$ is congruent to $1$ modulo $5$, the two congruences
in~(a) combine into $m^4\equiv5n-4\pmod{125}$. Whether $n$ lies in $\cN$ can be read
off from its prime factorization (Proposition~\ref{prop:g5-norms}), and whether a
prime $q$ lies in $\cB$ can be decided in $\F_q$ (Proposition~\ref{prop:g5-jacobi}),
so~(a) gives an effective test. Part~(b) shows that the bound of
Theorem~\ref{thm:A}(c) is attained for $p=5$ although $5^{13}$ itself is not a
determinant: valuation $13$ occurs only with a cofactor divisible by a prime in
$\cB$. For the abelian group of the same order and exponent, Mossinghoff and
Pinner~\cite[Theorem~3.2]{MossinghoffPinner2026} show that
$5^{22}\Z\subseteq S_0(C_5\times C_{25})\subseteq5^{13}\Z$ and that $5^{13}$ is not a
group determinant of $C_5\times C_{25}$.

\subsection*{Method}

The starting point is the Frobenius factorization (Lemma~\ref{lem:fact})
\[
 D_{\Gp}(F)=D_{C_p\times C_p}(\Fbar)\,\Nm_{k/\Q}\bigl(\Delta(F)\bigr)^p ,
\]
where $\Fbar$ is the image of $F$ in the abelianization $C_p\times C_p$,
$k=\Q(\omega)$, and $\Delta(F)\in\Z[\omega]$ is the determinant of $F$ in an
irreducible representation of degree $p$. The image of $\Z[\Gp]$ in that
representation is a crossed-product order $\Lam$ in $\M_p(k)$. On the congruence
subgroup $1+\pi\Lam_\p$, where $\pi=1-\omega$, its reduced norms fill exactly
$1+\pi^p\Z_p[\omega]$ (Lemma~\ref{lem:L1}), whereas for a maximal order they fill
$1+\pi\Z_p[\omega]$. With strong approximation for $\SL_p$, this shows that once
$\Fbar$ is fixed, $\Delta(F)$ can be any nonzero element allowed by one local
condition at $\pi$ (Proposition~\ref{prop:fibre}). The problem then becomes one of
deciding which values of the linear characters can occur together. Part~(a) of
Theorem~\ref{thm:A} comes from comparing two characters of
$\F_p[C_p\times C_p]^\times$ on an explicit basis. For $p=3$, global units and
integral Fourier inversion produce the required character values. For $p=5$ the cubic
invariants of the character values on the six lines through the origin of $\F_5^2$
must come from a single binary cubic form. This explains the exclusion of $418$, and,
together with a congruence for $\Delta(F)$ modulo $\pi^5$, the role of $\cB$ in
Theorem~\ref{thm:C}(b).

Section~\ref{sec:order} collects the representation theory and the lifting results.
Section~\ref{sec:U} proves parts~(a) and~(b) of Theorem~\ref{thm:A} together with the
norm-divisor condition, and Section~\ref{sec:B} proves part~(c) and
Corollary~\ref{cor:primes}. Sections~\ref{sec:p3}, \ref{sec:g5-coprime}
and~\ref{sec:g5-divisible} treat $p=3$, the coprime values for $p=5$, and the
multiples of $5$. Section~\ref{sec:open} ends with some questions for $p\ge7$.

We checked the explicit identities in the paper (the element $F_{57}$ of
Section~\ref{sec:U}, the elements of Section~\ref{sec:p3}, and the family
\eqref{eq:g5-explicit-family}) by exact computation in Magma and in Python: full
$27\times27$ and $125\times125$ determinants, and $\Delta(F_t)=\pi$ as an identity
in $t$. The proofs do not rely on
these computations. The scripts are available from the author.

\section{Representations and lifting}\label{sec:order}

Throughout, $p$ is an odd prime, $\zeta$ is a primitive $p^2$th root of unity,
$\omega=\zeta^p$, $K=\Q(\zeta)$ and $k=\Q(\omega)$. We write $\cO=\Z[\omega]$ and
$\cO_K=\Z[\zeta]$ for the rings of integers, and put $\pi=1-\omega$,
$\varpi=1-\zeta$, $\p=\pi\cO$ and $\PP=\varpi\cO_K$. These are the primes above $p$,
and $p\cO=\p^{p-1}$, $\pi\cO_K=\PP^p$. The automorphism $\sigma$ of $K$ with
$\sigma(\zeta)=\zeta^{1+p}=\zeta\omega$ generates $\Gal(K/k)$. For $a$ prime to $p$,
$\tau_a$ is the automorphism of $k$ with $\tau_a(\omega)=\omega^a$, and we fix an
extension of $\tau_a$ to $K$. We write $\Nm$ and $\Tr$ for the norm and trace from $k$
to $\Q$, and $\Nm_{K/k}$, $\Tr_{K/k}$ for the relative ones. Completions are written
with subscripts: $\cO_\p=\Z_p[\omega]$, $\cO_{\mathfrak q}$ is the completion at a
finite place $\mathfrak q$ of $k$, and $\cO_{K,\PP}$ is the completion of $\cO_K$ at
$\PP$.

Two facts about traces will be used repeatedly. Since
$\cO_K=\bigoplus_{t=0}^{p-1}\cO\zeta^t$ and $\Tr_{K/k}(\zeta^t)=0$ for $0<t<p$,
\begin{equation}\label{eq:trace-Kk}
 \Tr_{K/k}(\cO_K)=p\cO=\p^{p-1}.
\end{equation}
Since $\Tr(\omega^l)$ is $p-1$ if $p\mid l$ and $-1$ otherwise, an element
$\alpha=\sum_lc_l\omega^l$ with $c_l\in\Z$ has $\Tr(\alpha)=pc_0-\sum_lc_l$, and so
\begin{equation}\label{eq:trace-kQ}
 \Tr(\p)\subseteq p\Z,\qquad \Tr(p\p)=\Tr(\p^p)\subseteq p^2\Z,\qquad
 \p^{2p}\cap\Z=p^2\p^2\cap\Z\subseteq p^2\Z .
\end{equation}

\subsection{The group and its representations}

Every element of $\Gp$ can be written uniquely as $x^ay^b$ with $0\le a<p^2$ and
$0\le b<p$, and $yx^ay^{-1}=x^{a(1+p)}$. Since
$\sum_{l=0}^{p-1}(1+p)^{bl}\equiv p\pmod{p^2}$ for every $b$, we have
$(x^ay^b)^p=x^{ap}$, so $x^ay^b$ has order dividing $p$ exactly when $p\mid a$, and
these elements form the subgroup $N=\langle x^p,y\rangle\cong C_p\times C_p$. The
commutator $yxy^{-1}x^{-1}=x^p$ generates the derived subgroup, which is also the
centre $Z=\langle x^p\rangle$. We write
\[
 \Z[\Gp]\to\Z[\Gp/Z]=\Z[C_p\times C_p],\ F\mapsto\Fbar,\qquad
 \Z[\Gp]\to\Z[\Gp/N]=\Z[C_p],\ F\mapsto\Ftil,
\]
where $C_p\times C_p=\langle X\rangle\times\langle Y\rangle$ with $x\mapsto X$,
$y\mapsto Y$, and $\Gp/N=\langle X\rangle$ with $y\mapsto1$. The $p^2$ linear
characters of $\Gp$ are given by $\chi_{ij}(x)=\omega^i$ and $\chi_{ij}(y)=\omega^j$,
so that $\chi_{ij}(F)=\Fbar(\omega^i,\omega^j)$.

Let $V=K^p$ with basis $\mathbf e_0,\dots,\mathbf e_{p-1}$, indices read modulo $p$,
and put
\[
 \rho(x)\mathbf e_r=\sigma^r(\zeta)\,\mathbf e_r,\qquad
 \rho(y)\mathbf e_r=\mathbf e_{r-1}.
\]
Then $\rho(y)\rho(x)\rho(y)^{-1}\mathbf e_r=\sigma^{r+1}(\zeta)\mathbf e_r
=\rho(x)^{1+p}\mathbf e_r$, so $\rho$ is a representation of $\Gp$. It is induced
from the character $x\mapsto\zeta$ of $\langle x\rangle$, whose $p$ conjugates
$x\mapsto\sigma^r(\zeta)$ are distinct, and is therefore irreducible. Its character
vanishes off $\langle x\rangle$ and takes the value $\Tr_{K/k}(\zeta^a)$ at $x^a$;
this is $p\omega^{a/p}$ if $p\mid a$ and $0$ otherwise. Thus the character has values in
$k$, and the representations $\rho_a=\tau_a\circ\rho$, $1\le a\le p-1$, are the
irreducible representations of degree $p$; together with the $p^2$ linear characters
they account for $p^2+(p-1)p^2=p^3$. The matrix $u=\rho(y)$ has entries $1$ in the
positions $(r,r+1)$ and $0$ elsewhere, and for $F=\sum_{b=0}^{p-1}f_b(x)y^b$ with
$f_b\in\Z[x]/(x^{p^2}-1)$,
\begin{equation}\label{eq:rhoF}
 \rho(F)=\sum_b\diag\bigl(\sigma^rf_b(\zeta)\bigr)_ru^b,\qquad
 \rho(F)_{r,r+b}=\sigma^r\bigl(f_b(\zeta)\bigr).
\end{equation}
Thus $\rho(F)$ is a twisted circulant: its $b$th diagonal carries the
$\sigma$-conjugates of $f_b(\zeta)$. Applying $\sigma$ to every entry of $\rho(F)$
permutes the rows and columns cyclically, so $\Delta(F)=\det\rho(F)$ is fixed by
$\sigma$. Hence $\Delta(F)\in\cO$, and $\det\rho_a(F)=\tau_a(\Delta(F))$.

\begin{lemma}[Frobenius factorization]\label{lem:fact}
For every $F\in\Z[\Gp]$,
\[
 D_{\Gp}(F)=M_1(F)\,\Nm\bigl(\Delta(F)\bigr)^p,\qquad
 M_1(F):=D_{C_p\times C_p}(\Fbar)=\prod_{i,j=0}^{p-1}\Fbar(\omega^i,\omega^j).
\]
\end{lemma}

\begin{proof}
By Frobenius~\cite{Frobenius1896}, $D_G(F)=\prod_\rho\det\rho(F)^{\deg\rho}$, the
product taken over the irreducible representations of $G$. The linear characters
contribute $\prod_{i,j}\Fbar(\omega^i,\omega^j)$, which is the group determinant of
$\Fbar$ for $C_p\times C_p$ by Dedekind's theorem, and the representations of degree
$p$ contribute $\prod_a\tau_a(\Delta(F))^p=\Nm(\Delta(F))^p$.
\end{proof}

We split $M_1(F)$ according to whether the character is trivial on $y$:
\begin{equation}\label{eq:AP}
 A(F):=\prod_{i=0}^{p-1}\Fbar(\omega^i,1)=D_{C_p}(\Ftil),\qquad
 P(F):=\prod_{i=0}^{p-1}\prod_{j=1}^{p-1}\Fbar(\omega^i,\omega^j),
\end{equation}
so that $M_1(F)=A(F)P(F)$.
The equality $A(F)=D_{C_p}(\Ftil)$ holds because $\Fbar(X,1)=\Ftil(X)$. The Galois
orbits of the nontrivial characters $\chi_{ij}$ correspond to the $p+1$ lines through
the origin of $\F_p^2$. The characters with $j=0$ give $A(F)$, and $P(F)$ is the
product over the other $p$ lines of norms $\Nm(\beta)$ of elements $\beta\in\cO$; in
particular $P(F)$ is a positive integer when $M_1(F)\ne0$. If $p\nmid F(1)$, every
value $\Fbar(\omega^i,\omega^j)$ is congruent to $F(1)$ modulo $\p$, and so
$P(F)\equiv F(1)^{p(p-1)}\equiv1\pmod p$.

\subsection{The crossed-product order}

Let $\Lam=\rho(\Z[\Gp])\subset\M_p(\cO_K)$. By~\eqref{eq:rhoF},
$\Lam=\bigoplus_{b=0}^{p-1}\cO_Ku^b$, where $\beta\in\cO_K$ is identified with
$\diag(\sigma^r\beta)_r$ and $u\beta u^{-1}=\sigma(\beta)$. Thus $\Lam$ is the
crossed-product order of the cyclic extension $K/k$ with trivial cocycle. It is an
$\cO$-order in the central simple $k$-algebra $\mathcal A=\Lam\otimes\Q\cong\M_p(k)$,
whose reduced norm is the determinant, and $\nrd(\Lam)\subseteq\cO$. For a finite
place $\mathfrak q$ of $k$ we put $\Lam_{\mathfrak q}=\Lam\otimes_{\cO}\cO_{\mathfrak q}$.
The completion that matters most is $\Lam_\p$.

\begin{lemma}\label{lem:lifts}
Let $F\in\Z[\Gp]$. The elements $F'\in\Z[\Gp]$ with $\Fbar'=\Fbar$ are exactly the
elements $F'=F+(x^p-1)H$ with $H\in\Z[\Gp]$, and
$\{\rho(F'):\Fbar'=\Fbar\}=\rho(F)+\pi\Lam$.
\end{lemma}

\begin{proof}
The kernel of $\Z[\Gp]\to\Z[\Gp/Z]$ is the ideal generated by $x^p-1$, and
$\rho(x^p-1)=(\omega-1)I=-\pi I$, so
$\rho\bigl((x^p-1)\Z[\Gp]\bigr)=\pi\rho(\Z[\Gp])=\pi\Lam$.
\end{proof}

\begin{lemma}\label{lem:units}
Let $\mathfrak q$ be a finite place of $k$ and $a\in\Lam_{\mathfrak q}$. Then
$a\in\Lam_{\mathfrak q}^\times$ if and only if $\nrd(a)\in\cO_{\mathfrak q}^\times$.
\end{lemma}

\begin{proof}
If $a$ is a unit, so is $\nrd(a)$. Conversely, left multiplication by $a$ on the free
$\cO_{\mathfrak q}$-module $\Lam_{\mathfrak q}$ of rank $p^2$ has determinant
$\nrd(a)^p$, because on $\mathcal A\otimes k_{\mathfrak q}\cong\M_p(k_{\mathfrak q})$
it is the direct sum of $p$ copies of the action of $a$ on columns. If $\nrd(a)$ is a
unit, this map is bijective, so $ab=1$ for some $b\in\Lam_{\mathfrak q}$, and $a$ is a
unit of the finite $\cO_{\mathfrak q}$-algebra $\Lam_{\mathfrak q}$.
\end{proof}

\subsection{Local reduced norms}

For a square matrix $B$ let $e_\ell(B)$ be the $\ell$th elementary symmetric function
of its eigenvalues, so that $\det(I+tB)=\sum_\ell t^\ell e_\ell(B)$.

\begin{lemma}\label{lem:newton}
Let $B\in\Lam_\p$ and $s_m=\tr(B^m)$. Then for $1\le\ell\le p-1$,
\[
 e_\ell(B)\in\p^{p-1}\cO_\p=p\cO_\p\qquad\text{and}\qquad
 e_\ell(B)\equiv(-1)^{\ell-1}\frac{s_\ell}{\ell}\pmod{p\p},
\]
and $\nrd(1+\pi B)\in1+\p^p\cO_\p$. Moreover $\Nm(1+\delta)\equiv1\pmod{p^2}$ for
every $\delta\in\p^p\cO_\p$, so $\Nm\bigl(\nrd(1+\pi\Lam_\p)\bigr)\subseteq1+p^2\Z_p$.
The same statements hold for the commutative subring $\cO_{K,\PP}$ of $\Lam_\p$; in
particular $\Nm_{K/k}(1+\pi\beta)\in1+\p^p$ for $\beta\in\cO_{K,\PP}$.
\end{lemma}

\begin{proof}
For $B=\sum_b\beta_bu^b\in\Lam_\p$ the trace of $B$ is $\Tr_{K/k}(\beta_0)$, and
$B^m\in\Lam_\p$ for all $m$, so by~\eqref{eq:trace-Kk} every power sum
$s_m=\tr(B^m)$ lies in $p\cO_\p=\p^{p-1}$. The coefficients $e_\ell(B)$ of the
characteristic polynomial of $B$ lie in $\cO_\p$: they are integral, and they are
fixed by $\sigma$ because applying $\sigma$ to the entries conjugates $B$ by $u$.
Newton's identities $\ell e_\ell=\sum_{m=1}^\ell(-1)^{m-1}e_{\ell-m}s_m$, in which
$\ell$ is a unit for $\ell\le p-1$, give $e_\ell(B)\in\p^{p-1}$ by induction on
$\ell$. The terms with $1\le m<\ell$ are products of two elements of $p\cO_\p$, so
they lie in $p^2\cO_\p\subseteq p\p$, and the same identities give the stated
congruence. Therefore
$\det(1+\pi B)=1+\sum_{\ell=1}^{p-1}\pi^\ell e_\ell(B)+\pi^p\det B\in1+\p^p$.

For the norm, let $\delta\in\p^p$. Then $\Nm(1+\delta)=1+\Tr(\delta)+\sum_{m\ge2}e_m$,
where $e_m$ is the $m$th elementary symmetric function of the conjugates of $\delta$.
By the completed form of~\eqref{eq:trace-kQ}, $\Tr(\delta)\in p^2\Z_p$ and
$e_m\in\p^{mp}\cO_\p\cap\Z_p\subseteq p^2\Z_p$. The commutative case is the case
$B=\diag(\sigma^r\beta)$.
\end{proof}

For $\Lam_\p$ the bound of Lemma~\ref{lem:newton} is attained.

\begin{lemma}\label{lem:L1}
For every odd prime $p$,
\[
 \nrd(1+\pi\Lam_\p)=1+\p^p\cO_\p.
\]
\end{lemma}

\begin{proof}
Write $U^{(n)}=1+\pi^n\cO_\p$. The set $1+\pi\Lam_\p$ is a compact group, since
geometric series give inverses in the complete order, so its image $H$ under the
reduced norm is closed, and $H\subseteq U^{(p)}$ by Lemma~\ref{lem:newton}. As $u$ is
a cyclic permutation matrix of odd order $p$, the element $I+\pi u$ has determinant
$1+\pi^p$, so $H$ maps onto $U^{(p)}/U^{(p+1)}$. For $n>p$ put $j=n-p+1\ge2$. The
scalar matrix $(1+\pi^j)I$ has determinant $(1+\pi^j)^p$. In the binomial expansion
of $(1+\pi^j)^p-1$, the term $p\pi^j$ has valuation $p-1+j=n$, the terms with
$2\le r\le p-1$ have valuations $p-1+rj>n$, and the last term has valuation $pj>n$.
So $H$ meets $U^{(n)}\smallsetminus U^{(n+1)}$ for every $n\ge p$. Since each quotient
$U^{(n)}/U^{(n+1)}$ is cyclic of order $p$, any element of $U^{(p)}$ is the limit of
products of powers of these elements, and as $H$ is closed, $H=U^{(p)}$.
\end{proof}

\subsection{Lifting with a fixed abelianization}

The next proposition reduces the construction of group determinants to a local
question at $\p$.

\begin{proposition}\label{prop:fibre}
Let $F\in\Z[\Gp]$ and let $\delta\in\cO$ be nonzero. There is an $F'\in\Z[\Gp]$ with
$\Fbar'=\Fbar$ and $\Delta(F')=\delta$ if and only if
$\delta\in\nrd\bigl(\rho(F)+\pi\Lam_\p\bigr)$.
\end{proposition}

\begin{proof}
Necessity follows from Lemma~\ref{lem:lifts}. For sufficiency, fix an isomorphism
$\mathcal A\cong\M_p(k)$ and identify $\Lam$ with an $\cO$-order
$\Lam'\subset\M_p(k)$; the lattices $\Lam'$ and $\M_p(\cO)$ have the same completion
at all but finitely many places. Let $S$ be the finite set of finite places
$\mathfrak q$ of $k$ such that $\mathfrak q\mid p$, or $\mathfrak q\mid\delta$, or
$\Lam'_{\mathfrak q}\ne\M_p(\cO_{\mathfrak q})$. For each $\mathfrak q\in S$ we choose
$g_{\mathfrak q}\in\Lam'_{\mathfrak q}$ with $\nrd(g_{\mathfrak q})=\delta$. For
$\mathfrak q=\p$ this is possible by hypothesis. For $\mathfrak q\ne\p$ the extension
$K/k$ is unramified at $\mathfrak q$, because the discriminant of $\Q(\zeta)$ is a
power of $p$. So $\cO_{K,\mathfrak q}=\cO_K\otimes\cO_{\mathfrak q}$ is a product of
unramified extensions of $\cO_{\mathfrak q}$ with finite residue fields, and the norm
$\cO_{K,\mathfrak q}^\times\to\cO_{\mathfrak q}^\times$ is
surjective~\cite[Ch.~V, \S2]{Serre1979}. For $a,b\in\cO_{K,\mathfrak q}$ the element
$a-bu$ of $\Lam_{\mathfrak q}$ has determinant $\Nm_{K/k}(a)-\Nm_{K/k}(b)$: only the
identity and the $p$-cycle contribute to the expansion of the determinant, and the
$p$-cycle is even. (The identity is polynomial, so it holds in the \'etale algebra
$K\otimes k_{\mathfrak q}$.) Choosing $a$ with $\Nm_{K/k}(a)=1+\pi_{\mathfrak q}$ for a
uniformizer $\pi_{\mathfrak q}$ gives $\nrd(a-u)=\pi_{\mathfrak q}$. As $\nrd$ is
multiplicative and takes every unit value on the diagonal units, every nonzero element
of $\cO_{\mathfrak q}$ is a reduced norm from $\Lam_{\mathfrak q}$.

Put $g_0=\diag(\delta,1,\dots,1)\in\M_p(\cO)$, so that $\det g_0=\delta$, and for
$\mathfrak q\in S$ let
\[
 W_{\mathfrak q}=\{\gamma\in\SL_p(k_{\mathfrak q}):
 g_0\gamma\in g_{\mathfrak q}+\pi_{\mathfrak q}^{N}\Lam'_{\mathfrak q}\}
\]
with $N$ large. Each $W_{\mathfrak q}$ is open and contains
$g_0^{-1}g_{\mathfrak q}$. Since $\SL_p$ is $k$-simple and simply connected and
$\SL_p(\mathbf C)$ is not compact, strong approximation away from the archimedean
places~\cite[Theorem~7.12]{PlatonovRapinchuk1994} gives $\gamma\in\SL_p(k)$ with
$\gamma\in W_{\mathfrak q}$ for $\mathfrak q\in S$ and
$\gamma\in\SL_p(\cO_{\mathfrak q})$ for $\mathfrak q\notin S$. Then $g=g_0\gamma$
lies in $\Lam'_{\mathfrak q}$ for every finite $\mathfrak q$ (for $\mathfrak q\notin
S$ because $g_0\in\M_p(\cO_{\mathfrak q})=\Lam'_{\mathfrak q}$), so $g\in\Lam'=\Lam$,
and $\nrd(g)=\delta$. Taking $N\ge1$ at $\mathfrak q=\p$ gives
$g\in\rho(F)+\pi\Lam_\p$, and since $g$ and $\rho(F)$ both lie in $\Lam$, in fact
$g\in\rho(F)+\pi\Lam$. By Lemma~\ref{lem:lifts}, $g=\rho(F')$ for some $F'$ with
$\Fbar'=\Fbar$.
\end{proof}

\subsection{The sum of the group elements}

Let $\Sigma=\sum_{g\in\Gp}g$, and for $F\in\Z[\Gp]$ let
\[
 Q(F)=\prod_{(i,j)\ne(0,0)}\Fbar(\omega^i,\omega^j)\cdot\Nm\bigl(\Delta(F)\bigr)^p
\]
be the product of the factors in Lemma~\ref{lem:fact} that come from nontrivial
irreducible representations. Thus $D_{\Gp}(F)=F(1)\,Q(F)$.

\begin{lemma}\label{lem:shift}
For every $F\in\Z[\Gp]$ and $t\in\Z$,
\[
 D_{\Gp}(F+t\Sigma)=\bigl(F(1)+p^3t\bigr)\,Q(F).
\]
\end{lemma}

\begin{proof}
The element $\Sigma/p^3$ is the idempotent that projects onto the trivial
representation, so $\Sigma$ acts as $0$ in every nontrivial irreducible representation
and as $p^3$ in the trivial one, and $F+t\Sigma$ has the same nontrivial factors as
$F$.
\end{proof}

The same device was used for $H_3$ by Mossinghoff and
Pinner~\cite[proofs of Theorems~3.3 and~3.6]{MossinghoffPinnerHeisenberg}, and for
$C_5\times C_5$ in~\cite[Lemma~2.1]{Panraksa2026Order25}.

\subsection{The unit group of \texorpdfstring{$\F_p[C_p\times C_p]$}{Fp[Cp x Cp]}}

Let $R=\F_p[C_p\times C_p]=\F_p[X,Y]/(X^p-1,Y^p-1)$. With $e=1-X$ and $v=1-Y$ we have
$R=\F_p[e,v]/(e^p,v^p)$, a local ring whose maximal ideal $\m=(e,v)$ satisfies
$\m^{2p-1}=0$, and $R^\times=\F_p^\times\times(1+\m)$. For $F\in\Z[\Gp]$, every
character value is congruent to the augmentation $F(1)=\sum_ga_g$ modulo $\p$. So
$p\nmid D_{\Gp}(F)$ if and only if $p\nmid F(1)$, that is, if and only if
$\Fbar\bmod p$ lies in $R^\times$.

\begin{lemma}\label{lem:basis}
The group $1+\m$ is elementary abelian of order $p^{p^2-1}$, and the $p^2-1$ elements
$e_{ij}=1+(1-X)^i(1-Y)^j$, $0\le i,j\le p-1$, $(i,j)\ne(0,0)$, form a basis of it over
$\F_p$. Consequently two homomorphisms from $R^\times$ to an abelian group coincide if
they agree on $\F_p^\times$ and on the $e_{ij}$.
\end{lemma}

\begin{proof}
For $a=\sum c_{ij}e^iv^j\in\m$ we have $a^p=\sum c_{ij}^pe^{pi}v^{pj}=0$, so
$(1+a)^p=1$. For each $d\ge1$ the map $1+\m^d\to\m^d/\m^{d+1}$, $1+a\mapsto a$, is a
surjective homomorphism with kernel $1+\m^{d+1}$, because $(1+a)(1+b)=1+a+b+ab$ with
$ab\in\m^{2d}\subseteq\m^{d+1}$. The monomials $e^iv^j$ with $i+j=d$ and
$0\le i,j\le p-1$ form a basis of $\m^d/\m^{d+1}$, and they are the images of the
$e_{ij}$ with $i+j=d$. By descending induction on $d$, every element of $1+\m^d$ is a
product of powers of the $e_{ij}$ with $i+j\ge d$. Hence the $e_{ij}$ generate
$1+\m$, and $p^2-1$ generators of an elementary abelian group of rank
$p^2-1=\dim_{\F_p}\m$ form a basis.
\end{proof}

The same argument in one variable shows that
$\F_p[C_p]^\times=\F_p^\times\times\bigl(1+(1-X)\F_p[C_p]\bigr)$, where the second
factor is elementary abelian with basis $1+(1-X)^i$, $1\le i\le p-1$.

\section{The congruence modulo \texorpdfstring{$p^3$}{p\^{}3}}\label{sec:U}

In this section we prove parts~(a) and~(b) of Theorem~\ref{thm:A} and the
norm-divisor condition. Part~(a) is proved by comparing two homomorphisms from
$R^\times$ to $(\Z/p^3\Z)^\times$ on the basis of Lemma~\ref{lem:basis}. This needs
two computations, one for the basis elements $e_{i0}$, which lift to the cyclic
subgroup $\langle x\rangle$, and one for the others.

\subsection{Two computations}

On $\Z[\langle x\rangle]$ the group determinant of $\Gp$ is the $p$th power of a
circulant determinant of order $p^2$, so for the elements $e_{i0}$ we need a
congruence between circulant determinants of orders $p^2$ and $p$. Both computations
involve the integers
\[
 \kappa(n):=\sum_{p\mid l}(-1)^l\binom nl ,
\]
for which
\begin{equation}\label{eq:c}
 \Tr(\pi^n)=p\,\kappa(n),\qquad \Tr_{K/\Q}(\varpi^n)\equiv-p\,\kappa(n)\pmod{p^2}
 \qquad(n\ge1).
\end{equation} These follow from the values of $\Tr(\omega^l)$ given above and the fact
that $\Tr_{K/\Q}(\zeta^l)$ equals $p(p-1)$, $-p$ or $0$ according as $p^2\mid l$,
$p\mathrel\|l$ or $p\nmid l$. Note that $\kappa(n)$ is the constant coefficient of
$(1-X)^n$ in $\Z[C_p]$. For $1\le i\le p-1$ let $E_\ell(i)$, $1\le\ell\le p-1$, be
the elementary symmetric functions of the $p-1$ conjugates $(1-\omega^a)^i$, $a\ne0$,
of $\pi^i$. They are rational integers whose power sums are
$\Tr(\pi^{im})=p\,\kappa(im)$. Products of two power sums are divisible by $p^2$, so
Newton's identities modulo $p^2$ give
\begin{equation}\label{eq:E}
 E_\ell(i)\equiv(-1)^{\ell-1}\frac{p\,\kappa(i\ell)}{\ell}\pmod{p^2},\qquad
 \Nm(1+\pi^i)=1+\sum_{\ell=1}^{p-1}E_\ell(i)\equiv1+pH_i\pmod{p^2},
\end{equation}
where $\Nm(1+\pi^i)=\prod_{a=1}^{p-1}\bigl(1+(1-\omega^a)^i\bigr)$ and
$H_i=\sum_{\ell=1}^{p-1}(-1)^{\ell-1}\kappa(i\ell)/\ell\in\Z_p$.

\begin{lemma}\label{lem:I}
Let $g\in\Z[C_{p^2}]$ with $p\nmid g(1)$, and let $\bar g$ be its image in
$\Z[C_p]$. Then
\[
 D_{C_{p^2}}(g)\equiv D_{C_p}(\bar g)\pmod{p^2};\qquad\text{equivalently,}\qquad
 \Nm_{K/\Q}\bigl(g(\zeta)\bigr)\equiv1\pmod{p^2}.
\]
\end{lemma}

\begin{proof}
Write $g$ as a polynomial in a generator $x$ of $C_{p^2}$. Then
\[
 D_{C_{p^2}}(g)=\prod_{j=0}^{p^2-1}g(\zeta^j)
 =\prod_{i=0}^{p-1}\bar g(\omega^i)\cdot\prod_{p\nmid j}g(\zeta^j)
 =D_{C_p}(\bar g)\,\Theta(g),\qquad \Theta(g)=\Nm_{K/\Q}\bigl(g(\zeta)\bigr),
\]
so the two statements are equivalent, and $\Theta$ is multiplicative. We first show
that $\Theta(g)\bmod p^2$ depends only on the class of $g$ in $\F_p[C_p]$. Since
$g(1)$ is a unit, $g$ is a unit of the local ring $\Z_p[C_{p^2}]$, and an element
$g'$ with the same image in $\F_p[C_p]$ has the form
$g'=g\,(1+(x^p-1)h_1+ph_2)$ with $h_1,h_2\in\Z_p[C_{p^2}]$. Evaluating at $\zeta$,
and using $x^p-1\mapsto-\pi$ and $p\in\pi\cO_{K,\PP}$, we get
$g'(\zeta)=g(\zeta)(1-\pi\beta)$ with $\beta\in\cO_{K,\PP}$. By
Lemma~\ref{lem:newton}, $\Nm_{K/k}(1-\pi\beta)\in1+\p^p$, so
$\Theta(g')\equiv\Theta(g)\pmod{p^2}$. Hence $\Theta$ induces a character of
$\F_p[C_p]^\times$, and by the one-variable form of Lemma~\ref{lem:basis} it suffices
to show that $\Theta\equiv1\pmod{p^2}$ on constants $c\in\F_p^\times$ and on
$g=1+(1-x)^i$, $1\le i\le p-1$. For a constant, $\Theta(c)=c^{p(p-1)}\equiv1\pmod{p^2}$.

Fix $1\le i\le p-1$, put $a=\varpi^i=(1-\zeta)^i$, and write
$\Nm_{K/k}(1+a)=\prod_r(1+\sigma^ra)=1+\sum_{m=1}^{p}E'_m$, where $E'_m$ is the $m$th
elementary symmetric function of the $\sigma^ra$. Then
$E'_p=\Nm_{K/k}(a)=\pi^i$, because
$\prod_r(1-\sigma^r\zeta)=\prod_r(1-\zeta\omega^r)=1-\zeta^p=\pi$. The power sums
$s_m=\Tr_{K/k}(a^m)$ lie in $p\cO$ by~\eqref{eq:trace-Kk}. Newton's identities give
$E'_m\in p\cO$ for $m\le p-1$, and since the other terms in the identities are
products of two elements of $p\cO$,
\[
 E'_m\equiv(-1)^{m-1}\frac{s_m}{m}\pmod{p\p}\qquad(1\le m\le p-1).
\]
Let $\gamma'=\sum_{m<p}E'_m\in p\cO$. Then
$\Theta(1+(1-x)^i)=\Nm(1+\pi^i+\gamma')=\Nm(1+\pi^i)\,\Nm(1+\gamma'')$ with
$\gamma''=\gamma'/(1+\pi^i)\equiv\gamma'\pmod{p\p}$. Here $\gamma''\in p\cO_\p$, so
$\Nm(1+\gamma'')\equiv1+\Tr(\gamma'')\pmod{p^2}$, the higher elementary symmetric
functions of its conjugates lying in $p^2\Z_p$; and
$\Tr(\gamma'')\equiv\Tr(\gamma')\pmod{p^2}$ by the completed form
of~\eqref{eq:trace-kQ}. Therefore, by~\eqref{eq:c},
\[
 \Nm(1+\gamma'')\equiv1+\sum_{m=1}^{p-1}(-1)^{m-1}\frac{\Tr_{K/\Q}(\varpi^{im})}{m}
 \equiv1-pH_i\pmod{p^2},
\]
with $H_i$ as in~\eqref{eq:E}, while $\Nm(1+\pi^i)\equiv1+pH_i\pmod{p^2}$
by~\eqref{eq:E}. Hence
$\Theta(1+(1-x)^i)\equiv(1+pH_i)(1-pH_i)\equiv1\pmod{p^2}$.
\end{proof}

The elements $e_{ij}$ with $j\ge1$ lift to $F=1+(1-x)^i(1-y)^j$, and $\Ftil=1$. For
these the factors $P(F)$ and $\Nm(\Delta(F))^p$ in the Frobenius factorization are
both nontrivial modulo $p^3$, but they cancel.

\begin{lemma}\label{lem:II}
Let $1\le i,j\le p-1$, $F=1+(1-x)^i(1-y)^j\in\Z[\Gp]$, and
\[
 s=\sum_{\ell=1}^{p-1}(-1)^{\ell-1}\frac{\kappa(i\ell)\,\kappa(j\ell)}{\ell}\in\Z_p .
\]
Then $A(F)=1$, $P(F)\equiv1+p^2s\pmod{p^3}$ and
$\Nm(\Delta(F))\equiv1-ps\pmod{p^2}$. Consequently
$D_{\Gp}(F)\equiv1\pmod{p^3}$.
\end{lemma}

\begin{proof}
The character values are
$\Fbar(\omega^a,\omega^b)=1+(1-\omega^a)^i(1-\omega^b)^j$, so
$\Fbar(\omega^a,1)=1$ and $A(F)=1$. For $b\ne0$ the product over $a$ is the image under
$\omega\mapsto\omega^b$ of
\[
 \Pi=\prod_{a=0}^{p-1}\bigl(1+\pi^j(1-\omega^a)^i\bigr)
 =1+\sum_{\ell=1}^{p-1}\pi^{j\ell}E_\ell(i)=:1+\xi,
\]
with $E_\ell(i)\in p\Z$ as in~\eqref{eq:E} (the factor $a=0$ equals $1$). Hence
$P(F)=\Nm(\Pi)$, and $\xi\in p\p^j$. The elementary symmetric functions of order
$m\ge2$ of the conjugates of $\xi$ lie in $\p^{2(p-1+j)}\cap\Z\subseteq p^2\p\cap\Z=p^3\Z$,
so
\[
 P(F)\equiv1+\Tr(\xi)=1+\sum_\ell E_\ell(i)\Tr(\pi^{j\ell})\pmod{p^3}.
\]
With $E_\ell(i)\equiv(-1)^{\ell-1}p\,\kappa(i\ell)/\ell\pmod{p^2}$ and
$\Tr(\pi^{j\ell})=p\,\kappa(j\ell)$, this is $P(F)\equiv1+p^2s\pmod{p^3}$.

Now let $B=\rho\bigl((1-x)^i(1-y)^j\bigr)\in\Lam$, so that
$\Delta(F)=\det(I+B)=1+\sum_{\ell=1}^pe_\ell(B)$. The top coefficient vanishes:
$e_p(B)=\det B=\det\rho((1-x)^i)\det\rho((1-y)^j)$, and
$\det\rho(1-y)=\det(I-u)=0$. By Lemma~\ref{lem:newton}, $e_\ell(B)\in p\cO$ and
$e_\ell(B)\equiv(-1)^{\ell-1}\tr(B^\ell)/\ell\pmod{p\p}$ for $\ell\le p-1$. Put
$\gamma=\sum_{\ell<p}e_\ell(B)\in p\cO$. As in the previous proof,
$\Nm(\Delta(F))\equiv1+\Tr(\gamma)\pmod{p^2}$ and
\[
 \Tr(\gamma)\equiv\sum_{\ell=1}^{p-1}(-1)^{\ell-1}
 \frac{\Tr\bigl(\tr(B^\ell)\bigr)}{\ell}\pmod{p^2}.
\]
The number $\Tr(\tr(B^\ell))=\Tr\bigl(\tr\rho(b_\ell)\bigr)$, where
$b_\ell=((1-x)^i(1-y)^j)^\ell\in\Z[\Gp]$, is the value at $b_\ell$ of the character
$\psi=\sum_a\tr\rho_a$ of the sum of the representations of degree $p$. The regular
character of $\Gp$ is $\sum_{i,j}\chi_{ij}+p\psi$, and $\sum_{i,j}\chi_{ij}(g)$ is
$p^2$ for $g\in Z=\langle x^p\rangle$ and $0$ otherwise. So
$\psi(g)=p^2[g=1]-p[g\in Z]$; that is, $\psi$ takes the value $p(p-1)$ at the
identity, $-p$ on $Z\smallsetminus\{1\}$, and $0$ off $Z$. Therefore
\[
 \Tr\bigl(\tr(B^\ell)\bigr)=p^2(b_\ell)_1-p\sum_{z\in Z}(b_\ell)_z
 \equiv-p\,C_\ell\pmod{p^2},
\]
where $(b_\ell)_g$ is the coefficient of $g$ in $b_\ell$ and
$C_\ell=\sum_{z\in Z}(b_\ell)_z$ is the constant coefficient of the image
$\bar b_\ell=(1-X)^{i\ell}(1-Y)^{j\ell}$ of $b_\ell$ in $\Z[C_p\times C_p]$, namely
$\kappa(i\ell)\kappa(j\ell)$. Hence
$\Nm(\Delta(F))\equiv1-p\sum_\ell(-1)^{\ell-1}\kappa(i\ell)\kappa(j\ell)/\ell
=1-ps\pmod{p^2}$, and finally
\begin{align*}
 D_{\Gp}(F)=A(F)P(F)\Nm(\Delta(F))^p&\equiv(1+p^2s)(1-ps)^p\\
 &\equiv(1+p^2s)(1-p^2s)\equiv1\pmod{p^3}. \qedhere
\end{align*}
\end{proof}

\begin{remark}
By Lucas's theorem,
$\kappa(n)\equiv\sum_j(-1)^j\binom{\lfloor n/p\rfloor}{j}\pmod p$, which is $1$ if
$n<p$ and $0$ if $n\ge p$. So
$s\equiv\sum_{\ell=1}^{\nu}(-1)^{\ell-1}/\ell\pmod p$ with
$\nu=\lfloor(p-1)/\max(i,j)\rfloor$. For $p=3$ and $i=j=1$ this gives
$s\equiv1-\tfrac12\equiv2\pmod3$, and indeed
$D_{G_3}(1+(1-x)(1-y))=24470992=208\cdot49^3$, with
$P(F)=208\equiv1+9\cdot2\pmod{27}$ and $\Nm(\Delta(F))=49\equiv1-3\cdot2\pmod9$.
\end{remark}

\subsection{Proof of Theorem~\ref{thm:A}(a)}

Let $F\in\Z[\Gp]$ with $p\nmid F(1)$. We first show that $D_{\Gp}(F)\bmod p^3$
depends only on the class of $\Fbar$ in $R=\F_p[C_p\times C_p]$. As for every
$p$-group, $\F_p[\Gp]$ is a local ring with the augmentation as residue map; since
$p\nmid F(1)$, $F$ is a unit of $\Z_p[\Gp]$, and an element $F'\in\Z[\Gp]$ with the
same image in $R$ can be written $F'=F\,(1+pH+(1-x^p)H')$ with
$H,H'\in\Z_p[\Gp]$. Since $x^p$ is central,
$1+pH+(1-x^p)H'=(1+pH)(1+(1-x^p)H'')$ with $H''=(1+pH)^{-1}H'\in\Z_p[\Gp]$, and both
$D_{\Gp}$ and the factorization of Lemma~\ref{lem:fact} extend to $\Z_p[\Gp]$.

For the first factor, let $L_H$ be left multiplication by $H$ on $\Z_p[\Gp]$ and
$H_1$ the coefficient of the identity in $H$. Then
$D_{\Gp}(1+pH)=\sum_\ell p^\ell e_\ell(L_H)$ with $e_1(L_H)=\tr L_H=p^3H_1$ and
$e_2(L_H)=\frac12\bigl(\tr(L_H)^2-\tr(L_{H^2})\bigr)\in p^3\Z_p$, so
$D_{\Gp}(1+pH)\equiv1\pmod{p^3}$. For the second factor, the image of
$1+(1-x^p)H''$ in $\Z_p[C_p\times C_p]$ is $1$, so
\[
 D_{\Gp}\bigl(1+(1-x^p)H''\bigr)=\Nm\bigl(\nrd(I+\pi\rho(H''))\bigr)^p
 \in(1+p^2\Z_p)^p\subseteq1+p^3\Z_p
\]
by Lemma~\ref{lem:newton}. Hence $D_{\Gp}(F')\equiv D_{\Gp}(F)\pmod{p^3}$, and
$F\mapsto D_{\Gp}(F)\bmod p^3$ is a well-defined homomorphism
$\theta_1\colon R^\times\to(\Z/p^3\Z)^\times$. In the same way $A(F)=D_{C_p}(\Ftil)$
modulo $p^2$ depends only on $\Ftil\bmod p$, because
$D_{C_p}(1+pH)=\prod_a(1+pH(\omega^a))\equiv1\pmod{p^2}$. So
$F\mapsto A(F)^p\bmod p^3$ is a homomorphism
$\theta_2\colon R^\times\to(\Z/p^3\Z)^\times$, which factors through
$\Fbar\mapsto\Ftil$.

By Lemma~\ref{lem:basis} it suffices to compare $\theta_1$ and $\theta_2$ on
$\F_p^\times$ and on the $e_{ij}$. For a constant $c$,
$D_{\Gp}(c)=c^{p^3}\equiv c^{p^2}=A(c)^p\pmod{p^3}$. For $e_{ij}$ with $j\ge1$ we take
the lift $F=1+(1-x)^i(1-y)^j$, for which $\Ftil=1$ and $A(F)=1$. If $i\ge1$, then
$D_{\Gp}(F)\equiv1\pmod{p^3}$ by Lemma~\ref{lem:II}. If $i=0$, then
$F=1+(1-y)^j$ lies in $\Z[\langle y\rangle]$ and
$D_{\Gp}(F)=D_{C_p}(1+(1-y)^j)^{p^2}\in(1+p\Z)^{p^2}\subseteq1+p^3\Z$. Here we used
that left multiplication by an element of $\Z[J]$, for a subgroup $J\le G$, acts on
$\Z[G]=\bigoplus\Z[J]g$ as $[G:J]$ copies of its action on $\Z[J]$. For $e_{i0}$ with
$i\ge1$ we take $F=1+(1-x)^i\in\Z[\langle x\rangle]$, so that
$D_{\Gp}(F)=D_{C_{p^2}}(F)^p$, while $A(F)=D_{C_p}(\Ftil)$ with $\Ftil$ the image of
$F$ in $\Z[C_p]$. Lemma~\ref{lem:I} gives
$D_{C_{p^2}}(F)\equiv D_{C_p}(\Ftil)\pmod{p^2}$, and hence
$D_{\Gp}(F)\equiv A(F)^p\pmod{p^3}$. Thus $\theta_1=\theta_2$, which is
Theorem~\ref{thm:A}(a). \qed

\begin{remark}\label{rem:siblings}
The proof also shows that, for every odd $p$, a group determinant of $\Gp$ prime to
$p$ is congruent modulo $p^3$ to $D_{C_{p^2}}(F(x,1))^p$, where $F(x,1)\in\Z[C_{p^2}]$
is obtained by substituting $y=1$ (which is not a ring homomorphism). Indeed
$D_{C_{p^2}}(F(x,1))\equiv D_{C_p}(\Ftil)\pmod{p^2}$ by Lemma~\ref{lem:I}, since
$F(x,1)$ reduces to $\Ftil$ modulo $x^p-1$.

For the abelian group $C_p\times C_{p^2}$, write $F(u,v)\in\Z[C_p\times C_{p^2}]$ with
$u^p=v^{p^2}=1$ and $p\nmid F(1,1)$, and put $B=D_{C_{p^2}}(F(1,v))$. Theorem~2.1 of
De Silva, Mossinghoff, Pigno and Pinner~\cite{DSMPP2019} (with $\ell=1$ and $k=2$ in their notation)
gives $D_{C_p\times C_{p^2}}(F)\equiv B^p\pmod{p^3}$, and Lemma~\ref{lem:I} gives
$B\equiv D_{C_p}(\Ftil)\pmod{p^2}$, where $\Ftil$ is the image of $F$ in the quotient
by the elements of order dividing $p$. Raising to the $p$th power gives the analogue
of Theorem~\ref{thm:A}(a) for $C_p\times C_{p^2}$. Its case $p=3$ is used in the
proof of Theorem~3.4 of~\cite{MossinghoffPinner2026}.

For the Heisenberg group $H_p$ the corresponding quotient is trivial. The scalar
element $F=(1+p)\cdot1$ has $D_{H_p}(F)=(1+p)^{p^3}\equiv1\pmod{p^3}$, whereas
$D_{\{1\}}(\Ftil)^p=(1+p)^p\equiv1+p^2\pmod{p^3}$, so the analogous congruence fails.
\end{remark}

\subsection{The norm divisor}

The necessary condition needs one more relation between the two factors of
$M_1(F)=A(F)P(F)$.

\begin{lemma}\label{lem:P}
Let $F\in\Z[\Gp]$ with $p\nmid F(1)$. Then $P(F)\equiv A(F)^{p-1}\pmod{p^2}$, and
consequently $M_1(F)\equiv A(F)^p\pmod{p^2}$.
\end{lemma}

\begin{proof}
Both sides are characters of $R^\times$. For $P$: if $F'=F(1+pH)$ with
$H\in\Z_p[\Gp]$, then
$P(F')=P(F)\prod_{j\ne0}\prod_i\bigl(1+p\,\overline H(\omega^i,\omega^j)\bigr)$ and
$\sum_{j\ne0}\sum_i\overline H(\omega^i,\omega^j)
=p^2\overline H_1-p\sum_b\overline H_{(0,b)}\in p\Z_p$, so
$P(F')\equiv P(F)\pmod{p^2}$; and $P$ is unchanged by $F\mapsto F+(1-x^p)H'$. For
$A^{p-1}$ this was shown in the proof of Theorem~\ref{thm:A}(a). On constants both
sides are congruent to $1$. On $e_{i0}$ the values
$\Fbar(\omega^a,\omega^b)=1+(1-\omega^a)^i$ do not depend on $b$, so
$P(F)=A(F)^{p-1}$ exactly. On $e_{ij}$ with $i,j\ge1$ we have $A(F)=1$ and
$P(F)\equiv1\pmod{p^2}$ by Lemma~\ref{lem:II}, and on $e_{0j}$ we have $A(F)=1$ and
$P(F)=\Nm(1+\pi^j)^p\equiv1\pmod{p^2}$ because $\Nm(1+\pi^j)\equiv1\pmod p$.
\end{proof}

\begin{proposition}\label{prop:nec}
Let $p$ be an odd prime and let $m\in S(\Gp)$ with $p\nmid m$. Then $m\in T_2$, and
$m$ has a norm divisor $n$ with $mn^p\in T_3$.
\end{proposition}

\begin{proof}
Let $m=D_{\Gp}(F)$, and write $A=A(F)$ and $P=P(F)$, so that
$m=AP\,\Nm(\Delta(F))^p$ by Lemma~\ref{lem:fact} and~\eqref{eq:AP}. The group
$(\Z/p^3\Z)^\times$ is the product of a cyclic group of order $p-1$ and one of order
$p^2$. Its subgroup of $p$th powers has index $p$ and is contained in the subgroup of
index $p$ formed by the classes of $T_2$, so the two coincide, and
Theorem~\ref{thm:A}(a) gives $m\in T_2$. Let $n=P$. By the remarks
after~\eqref{eq:AP}, $n$ is a positive integer and a product of norms from
$\Z[\omega]$, hence itself a norm from $\Z[\omega]$; and $n$ divides $M_1(F)$, which
divides $m$. By Lemma~\ref{lem:P}, $n^p\equiv A^{p(p-1)}\pmod{p^3}$, so by
Theorem~\ref{thm:A}(a)
\[
 mn^p\equiv A^p\cdot A^{p(p-1)}=A^{p^2}\pmod{p^3}.
\]
Since $A^{p^2}$ is the Teichm\"uller representative of $A$ modulo $p^3$, we get
$mn^p\in T_3$.
\end{proof}

\subsection{Proof of Theorem~\ref{thm:A}(b)}

Put $\Phi_p(z)=1+z+\cdots+z^{p-1}$. By Lemma~\ref{lem:shift} it is enough to find, for
each of the $p-1$ residue classes of $T_3$ modulo $p^3$, an element $F$ in that class
with $Q(F)=1$. We obtain $\Delta(F)$ as the determinant of a circulant perturbed by a
matrix of rank one; the coefficients of the element $h$ below are chosen to control
both the abelianization and the relative trace that appears in this perturbation.

Fix $a\in\{1,\dots,p-1\}$ and put $f_a(z)=1+z+\cdots+z^{a-1}$. In $\Z[X]/(X^p-1)$
define $q(X)=\sum_{i=0}^{p-1}q_iX^i$ by
\[
 p\,q(X)=f_a(X)^p-a .
\]
The coefficients $q_i$ are integers, because $f_a(X)^p\equiv\sum_{j=0}^{a-1}X^{jp}=a$
modulo $p$. Put $s=q(1)-q_0$ and
\[
 h(x)=q(x^p)+\sum_{i=1}^{p-1}(q_i+s)x^i,\qquad
 F_a=f_a(y)+h(x)\Phi_p(y)\in\Z[\Gp].
\]
In the abelianization $x^p$ becomes $1$, and the constant coefficient of $\bar h$ is
$q(1)=q_0+s$. Thus
\[
 \bar h(X)=q(X)+s\Phi_p(X),\qquad F_a(1)=a+p\bigl(q(1)+ps\bigr)=a^p+p^2s=:M_a .
\]
For a nontrivial linear character $(i,j)\ne(0,0)$,
\[
 \overline F_a(\omega^i,\omega^j)=
 \begin{cases}
 f_a(\omega^j),&j\ne0,\\
 f_a(\omega^i)^p,&j=0,\ i\ne0.
 \end{cases}
\]
Each $f_a(\omega^r)=(1-\omega^{ar})/(1-\omega^r)$ is an algebraic integer of norm $1$,
because multiplication by $a$ permutes the nonzero residues modulo $p$. Thus all the
nontrivial linear factors have norm $1$.

For the factor of degree $p$, let $\mathbf1$ be the column vector of $p$ ones, so
that $\mathbf1\mathbf1^t=\Phi_p(u)$, and put $v=(h(\sigma^r\zeta))_{0\le r<p}$. Then
\[
 \rho(F_a)=f_a(u)+v\mathbf1^t .
\]
The circulant $B=f_a(u)$ has eigenvalues $a$ and the $p-1$ conjugates of
$f_a(\omega)$, so $\det B=a$; moreover $\mathbf1^tB=a\mathbf1^t$. The determinant
formula for a rank-one perturbation gives
\[
 \Delta(F_a)=\det B\bigl(1+\mathbf1^tB^{-1}v\bigr)=a+\Tr_{K/k}\bigl(h(\zeta)\bigr).
\]
The terms $\zeta^i$, $1\le i<p$, have trace zero, while $q(\zeta^p)=q(\omega)\in k$.
Hence
\[
 \Delta(F_a)=a+p\,q(\omega)=f_a(\omega)^p ,
\]
which again has norm $1$. Therefore $Q(F_a)=1$ and $D_{\Gp}(F_a)=M_a$; the character
calculation with $j=0$ also gives $A(F_a)=M_a$. Since $M_a\equiv a\pmod p$,
Theorem~\ref{thm:A}(a) gives $M_a\equiv M_a^p\pmod{p^3}$, so $M_a\in T_3$.

For each nonzero residue $a$ modulo $p$, the equation $z^{p-1}=1$ has exactly one
solution modulo $p^3$ congruent to $a$, by Hensel's lemma. So the $M_a$ represent
the $p-1$ classes of $T_3$, and by Lemma~\ref{lem:shift} the progressions
$D_{\Gp}(F_a+t\Sigma)=M_a+p^3t$, $t\in\Z$, contain all of $T_3$. \qed

\medskip
For example, when $p=5$ the element
\[
 F_{57}=y^2+y^3+\bigl(-1+2x+3x^2+3x^3+2x^4+x^{10}+x^{15}\bigr)\Phi_5(y)
\]
has augmentation $57$. With $b=\omega^2+\omega^3$, its nontrivial linear values are
the conjugates of $b$ when $j\ne0$ and of $b^5$ when $j=0$, all of norm $1$, and the
same rank-one calculation gives
\[
 \Delta(F_{57})=2+5(-1+\omega^2+\omega^3)=-3+5\omega^2+5\omega^3=(1+\omega)^5 .
\]
By Lemma~\ref{lem:shift}, $D_{G_5}(F_{57}+t\Sigma)=57+125t$. Together with $1+t\Sigma$ and the negatives,
this realizes the four classes $1,57,68,124$ modulo $125$ directly.

\section{Divisibility and prime values}\label{sec:B}

\subsection{Proof of Theorem~\ref{thm:A}(c)}

The lower bound comes from the two factors in Lemma~\ref{lem:fact} separately, and the
sharpness from a first-order computation of $\Delta$ for an explicit element.

Suppose that $p\mid D_{\Gp}(F)$. Modulo $\PP$ every entry of $\rho(F)$ is congruent to
a rational integer, and the matrix becomes a circulant over $\F_p$ whose eigenvalues
all reduce to the row sum; so $\Delta(F)\equiv F(1)^p\pmod\p$. Every linear character
value is congruent to $F(1)$ modulo $\p$ as well. Hence
$D_{\Gp}(F)\equiv F(1)^{p^3}\pmod p$, and $p\mid D_{\Gp}(F)$ means
$p\mid F(1)=\Fbar(1,1)$. Then $p^{p+3}\mid M_1(F)=D_{C_p\times C_p}(\Fbar)$ by
Theorem~2.1 of~\cite{MossinghoffPinnerHeisenberg}, which holds for every odd $p$ (this
bound already appears in Trafford~\cite[Ch.~2, Lemma~13]{Trafford1992}), and
$\p\mid\Delta(F)$ gives $p^p\mid\Nm(\Delta(F))^p$. Hence $p^{2p+3}\mid D_{\Gp}(F)$.

For sharpness when $p\ge5$, let $F=p+(1-x)+(1-y)^2$, so that
$\Fbar=p+(1-X)+(1-Y)^2$. Then $\Fbar(1,1)=p$. For $i\ne0$ the value
$\Fbar(\omega^i,\omega^j)$ has $\p$-adic valuation exactly $1$, since $1-\omega^i$
dominates; and for $i=0$, $j\ne0$ the value $p+(1-\omega^j)^2$ has valuation exactly
$2$, because $p\in\p^{p-1}\subseteq\p^4$. Taking norms, the trivial character contributes $1$
to the $p$-adic valuation of $M_1(F)$, each of the $p$ orbits with $i\ne0$
contributes $1$, and the orbit with $i=0$ contributes $2$; so $M_1(F)$ has $p$-adic
valuation $p+3$, as in the proof of Theorem~2.1 of~\cite{MossinghoffPinnerHeisenberg}. It remains to show that $\p$
divides $\Delta(F)$ exactly once. Since $\sigma^r(\zeta)=\zeta\omega^r$ and
$\omega^r=(1-\pi)^r\equiv1-r\pi\pmod{\pi^2}$, we have
$1-\sigma^r(\zeta)\equiv\varpi+r\zeta\pi\pmod{\pi^2\cO_K}$, so that
\[
 \rho(F)\equiv C+\pi\zeta D\pmod{\pi^2\M_p(\cO_K)},
\]
where $C=(p+\varpi)I+(I-u)^2$ and $D=\diag(0,1,\dots,p-1)$.
Expanding the determinant to first order in the perturbation, whose entries lie in
$\pi\cO_K$, we obtain
\[
 \Delta(F)\equiv\det C+\pi\,\tr\bigl(\adj(C)\,\zeta D\bigr)\pmod{\pi^2\cO_K}.
\]
The matrix $C$ is a polynomial in $u$, hence so is $\adj(C)$, and a polynomial in $u$
has constant diagonal. If $\alpha\in\cO_K$ is that diagonal entry, then
$\tr(\adj(C)\zeta D)=\zeta\alpha\sum_{r<p}r=\zeta\alpha\,p(p-1)/2\in p\cO_K
\subseteq\pi\cO_K$. The eigenvalues of the circulant $C$ are
$p+\varpi+(1-\omega^j)^2$, $0\le j<p$, and each is congruent to $\varpi$ modulo
$\PP^{2p}=\pi^2\cO_K$, since $p\in\PP^{p(p-1)}$ and $(1-\omega^j)^2\in\PP^{2p}$.
Hence $\det C\equiv\varpi^p\equiv1-\zeta^p=\pi\pmod{\pi^2\cO_K}$, the second
congruence because $(1-\zeta)^p\equiv1-\zeta^p\pmod{p\cO_K}$. Therefore
$\Delta(F)\equiv\pi\pmod{\pi^2\cO_K}$, so $\Delta(F)\in\cO$ has $\p$-adic valuation
$1$, and $p^{2p+3}$ exactly divides $D_{\Gp}(F)$.

For $p=3$ a direct computation gives $D_{G_3}(x^4+x^7+x^6y^2)=3^9$, with $M_1=3^6$
and $\Nm(\Delta)=3$. (Alternatively, the element $F_1$ of Section~\ref{sec:p3-mult}
has determinant $3^9$ by~\eqref{eq:G3-progressions} with $t=0$.) \qed

\begin{remark}\label{rem:MP-exponent}
The proof follows the pattern of Theorem~3.2(i) of~\cite{MossinghoffPinner2026}
for $C_p\times C_{p^2}$, where the exponent is also $2p+3$: here $p+3$ comes from the
elementary abelian quotient and the remaining $p$ from the exponent $p$ on
$\Nm(\Delta(F))$. For $p\ge5$, part~(ii) of that theorem shows that $p^{2p+3}$ itself
is not a group determinant of $C_p\times C_{p^2}$, while for $p=3$ both groups attain
$3^9$.
\end{remark}

\subsection{Proof of Corollary~\ref{cor:primes}}

A prime $q\in S(\Gp)$ is different from $p$ by Theorem~\ref{thm:A}(c), so
Proposition~\ref{prop:nec} applies with $n\in\{1,q\}$. If $n=1$, then $q\in T_3$. If
$n=q$, then $q\in T_2$ and $q^{p+1}\in T_3$. Write $q\equiv\hat q(1+p^2t)\pmod{p^3}$,
where $\hat q=q^{p^2}$ is the Teichm\"uller representative of $q$ modulo $p^3$. Then
$q^{p+1}\equiv\hat q^{\,p+1}(1+p^2t)\pmod{p^3}$, and $q^{p+1}\in T_3$ forces $p\mid t$,
that is, $q\in T_3$. Conversely, Theorem~\ref{thm:A}(b) realizes every prime in $T_3$.

For the Lind--Lehmer constant, $D_{\Gp}(1+x)=D_{C_{p^2}}(1+x)^p=2^p$, since the
circulant determinant of $1+x$ of odd order $n$ is $\prod_{\eta^n=1}(1+\eta)=2$. For
$p=3$, every element of $S^*(G_3)$ other than $\pm1$ has absolute value at least $8$, because
$S^*(G_3)\subseteq T_2=\{m:m\equiv\pm1\pmod9\}$, and every nonzero multiple of $3$ in
$S(G_3)$ has absolute value at least $3^9$ by Theorem~\ref{thm:A}(c); hence
$\lambda(G_3)=\frac1{27}\log8$. For $p=5$, the elements of $T_2$ in $[2,31]$ are
$7,18,24,26$. None of them is in $T_3$, and none has a norm divisor $n>1$, because
each of the primes $2,3,7,13$ has order $4$ modulo $5$ and would have to occur to a
power divisible by $4$. So none of them lies in $S(G_5)$ by
Proposition~\ref{prop:nec}, and the nonzero multiples of $5$ in $S(G_5)$ have
absolute value at least $5^{13}$ by Theorem~\ref{thm:A}(c). Since
$D_{G_5}(-F)=-D_{G_5}(F)$, the same holds on $[-31,-2]$, and
$\lambda(G_5)=\frac1{125}\log32$. \qed

\section{The group of order \texorpdfstring{$27$}{27}}\label{sec:p3}

In this section, $p=3$, so $\omega$ is a primitive cube root of unity,
$\cO=\Z[\omega]$, $\pi=1-\omega$ and $\cO_\p=\Z_3[\omega]$.

\subsection{Coprime values}

Lemma~\ref{lem:L1} and Proposition~\ref{prop:fibre} reduce the `if' direction of
Theorem~\ref{thm:B}(a) to two facts special to $p=3$: a local unit whose norm is $1$
modulo $9$ is congruent to a global unit modulo $\p^3$, and $M_1$ can be prescribed
together with its factorization $A\cdot P$.

\begin{lemma}\label{lem:kernel}
The norm induces a homomorphism $(\cO/\p^3)^\times\to(\Z/9\Z)^\times$ with image
$\{1,4,7\}$, whose kernel consists of the classes of the six units
$\pm1,\pm\omega,\pm\omega^2$. In particular, if $\alpha\in\cO$ is prime to $3$ and
$\Nm(\alpha)\equiv1\pmod9$, then $\alpha\equiv\eta\pmod{\p^3}$ for a unit $\eta$ of
$\cO$.
\end{lemma}

\begin{proof}
The ring $\cO/\p^3$ has $27$ elements, $18$ of which are units. The map is well
defined because $\Nm(1+\p^3)\subseteq1+9\Z$ by
Lemma~\ref{lem:newton}. Norms of elements prime to $3$ are congruent to $1$ modulo
$3$, and $\Nm(1)=1$, $\Nm(1-3\omega)=13\equiv4$, $\Nm(2-\omega)=7$, so the image is
$\{1,4,7\}$ and the kernel has order $18/3=6$. The six units have norm $1$ and are
pairwise incongruent modulo $3=\p^2$, so they make up the kernel.
\end{proof}

\begin{lemma}\label{lem:atoms}
Let $m\in\Z$ with $m\equiv\eta\pmod9$, $\eta\in\{\pm1\}$, and let $n$ be a norm
divisor of $m$. Then there is an $\Fbar\in\Z[C_3\times C_3]$ with
\[
 \Fbar(1,1)=m/n,\qquad \Fbar(1,\omega)=\gamma,\qquad
 \Fbar(\omega,1)=\Fbar(\omega,\omega)=\Fbar(\omega,\omega^2)=\eta
\]
for some $\gamma\in\cO$ with $\Nm(\gamma)=n$. Every $F\in\Z[G_3]$ lifting $\Fbar$
has $M_1(F)=m$, $A(F)=m/n$ and $P(F)=n$.
\end{lemma}

\begin{proof}
Choose $\gamma_0\in\cO$ with $\Nm(\gamma_0)=n$. Since $(\cO/3\cO)^\times$ has order
$6$ and the six units of $\cO$ are pairwise incongruent modulo $3$, we may multiply
$\gamma_0$ by a unit to obtain $\gamma=r+s\omega\equiv\eta\pmod3$, that is,
$r\equiv\eta$ and $s\equiv0\pmod3$. Put $c_0=m/n$. The prescribed values, completed
by complex conjugation ($\Fbar(1,\omega^2)=\bar\gamma$, the three values $\eta$
being real), determine $\Fbar=\sum_{a,b}f_{ab}X^aY^b$ through the inverse Fourier
transform $f_{ab}=\frac19\sum_{i,j}\Fbar(\omega^i,\omega^j)\omega^{-(ia+jb)}$, and a
direct evaluation gives
\[
 f_{a0}=\frac{c_0+2r-s+3\eta(3\delta_{a0}-1)}9,\qquad
 f_{a1}=\frac{c_0-r+2s}9,\qquad f_{a2}=\frac{c_0-r-s}9 ,
\]
where $\delta_{a0}$ is $1$ for $a=0$ and $0$ otherwise. These nine numbers are
integers if and only if $3\mid s$, $r\equiv\eta\pmod3$ and $c_0\equiv r+s\pmod9$. The
first two conditions hold by the choice of $\gamma$. For the third, note that
$(r+s)\Nm(\gamma)=(r+s)(r^2-rs+s^2)=r^3+s^3\equiv\eta\pmod9$ and
$c_0\Nm(\gamma)=m\equiv\eta\pmod9$; since $\Nm(\gamma)\equiv1\pmod3$ is invertible
modulo $9$, $c_0\equiv r+s\pmod9$. Finally
$A(F)=\Fbar(1,1)\Fbar(\omega,1)\Fbar(\omega^2,1)=c_0\eta^2=m/n$ and
$P(F)=\Nm(\gamma)\eta^4=n$, so $M_1(F)=m$.
\end{proof}

\begin{example}
Take $m=62\equiv-1\pmod9$ and $n=31=\Nm(5+6\omega)$. Then $\eta=-1$,
$\gamma=5+6\omega\equiv-1\pmod3$ and $c_0=2$, and the lemma produces
$\Fbar=X+X^2+(1+X+X^2)(Y-Y^2)$. One lift of $\Fbar$ with $\Delta=-1$, and hence with
$D_{G_3}(F)=62$, is
\begin{multline*}
 F=-y^2+2x+2xy-xy^2-x^2-3x^2y^2+2x^3+x^3y-x^3y^2-2x^4y\\
 +2x^5-x^5y-x^5y^2-2x^6+x^6y^2-x^7+x^7y+2x^8y+3x^8y^2 .
\end{multline*}
(This lift was found by a computer search.)
\end{example}

\begin{proof}[Proof of Theorem~\ref{thm:B}(a)]
Let $m$ be an integer with $3\nmid m$. If $m\in S(G_3)$, then
Proposition~\ref{prop:nec} gives $m\in T_2=\{m:m\equiv\pm1\pmod9\}$ and a norm
divisor $n$ of $m$ with $mn^3\in T_3=\{m:m\equiv\pm1\pmod{27}\}$.

Conversely, suppose that $m\equiv\eta\pmod9$ with $\eta=\pm1$, and that $n$ is a norm
divisor of $m$ with $mn^3\equiv\pm1\pmod{27}$. Since $n\equiv1\pmod3$, the sign is
$\eta$, so $mn^3\equiv\eta\pmod{27}$. Lemma~\ref{lem:atoms} gives
$\Fbar\in\Z[C_3\times C_3]$ with $M_1=m$, $A=m/n$ and $P=n$. Let $F'\in\Z[G_3]$ be a
lift and $\Delta_0=\Delta(F')$, so that $D_{G_3}(F')=m\Nm(\Delta_0)^3$ and
$\Delta_0$ is prime to $3$. By Theorem~\ref{thm:A}(a),
\[
 m\Nm(\Delta_0)^3\equiv\Bigl(\frac mn\Bigr)^3\pmod{27},\qquad\text{so}\qquad
 \Nm(\Delta_0)^3\equiv m^2n^{-3}\equiv\eta m^3\equiv1\pmod{27},
\]
using $n^{-3}\equiv\eta m$ and $m^3\equiv\eta\pmod{27}$. As
$\Nm(\Delta_0)\equiv1\pmod3$, we may write $\Nm(\Delta_0)=1+3t$; then
$\Nm(\Delta_0)^3\equiv1+9t\pmod{27}$, so $3\mid t$, that is,
$\Nm(\Delta_0)\equiv1\pmod9$. Lemma~\ref{lem:kernel} now gives a unit $\eta_0$ of
$\cO$ with $\Delta_0\equiv\eta_0\pmod{\p^3}$. Then
$\eta_0/\Delta_0\in1+\p^3\cO_\p$, so by Lemma~\ref{lem:L1} there is
$W\in1+\pi\Lam_\p$ with $\nrd(W)=\eta_0/\Delta_0$, and
$\rho(F')W\in\rho(F')+\pi\Lam_\p$ has reduced norm $\eta_0$. By
Proposition~\ref{prop:fibre} there is a lift $F$ of $\Fbar$ with
$\Delta(F)=\eta_0$, and $D_{G_3}(F)=M_1(F)\Nm(\eta_0)^3=m$.
\end{proof}

\subsection{The multiples of \texorpdfstring{$3$}{3}}\label{sec:p3-mult}

Now $\zeta$ is a primitive ninth root of unity and $\omega=\zeta^3$. By
Theorem~\ref{thm:A}(c) and Lemma~\ref{lem:shift}, Theorem~\ref{thm:B}(b) will follow
once we find elements $F_c$ with $F_c(1)=3c$ for $c=0,1,2,3,4$ and $Q(F_c)=3^8$.
Consider
\begin{align*}
 F_0&=-x^3+(1-x+x^2)y^2,\\
 F_1&=x+x^2-y+(1+x)y^2,\\
 F_2&=x^2+(1+2x+x^2)y+x^2y^2,\\
 F_3&=1+x+2x^2+x^2y+(1+2x+x^2)y^2,\\
 F_4&=2+x+2x^2+(1+x+x^2-x^7)y+(2+2x+x^2)y^2.
\end{align*}
Their augmentations are $F_c(1)=3c$. The table lists their remaining factors, with
one column for each of the four nontrivial Galois orbits of linear characters.
\[
\begin{array}{c rrrr r}\toprule
 c&\overline F_c(1,\omega)&\overline F_c(\omega,1)&
 \overline F_c(\omega,\omega)&\overline F_c(\omega,\omega^2)&\Delta(F_c)\\ \midrule
 0&-2-\omega&-1-2\omega&-3&1+2\omega&\omega-1\\
 1&-3\omega&\omega-1&-1-2\omega&\omega-1&2+\omega\\
 2&3\omega&-2-\omega&-2-\omega&1-\omega&1+2\omega\\
 3&-3\omega&-2-\omega&1-\omega&-2-\omega&1-\omega\\
 4&-3\omega&1-\omega&1-\omega&-2-\omega&1+2\omega\\ \bottomrule
\end{array}
\]
The character values follow by substitution. The last column can be checked without
leaving $k$, as follows. Write $F_c=f_0(x)+f_1(x)y+f_2(x)y^2$, put $z_b=f_b(\zeta)$,
and recall $\sigma(\zeta)=\zeta^4$. Expanding the $3\times3$ determinant
in~\eqref{eq:rhoF} gives
\begin{equation}\label{eq:cubic-S0}
 \Delta(F_c)=\Nm_{K/k}(z_0)+\Nm_{K/k}(z_1)+\Nm_{K/k}(z_2)
 -\Tr_{K/k}\bigl(z_0\sigma(z_1)\sigma^2(z_2)\bigr),
\end{equation}
and for $A,B,C\in k$ we have $\Tr_{K/k}(A+B\zeta+C\zeta^2)=3A$ and
\[
 \Nm_{K/k}(A+B\zeta+C\zeta^2)=A^3+\omega B^3+\omega^2C^3-3\omega ABC .
\]
For example, for $F_4$ the three norms are $-19\omega$, $-3-7\omega$ and $7-5\omega$,
and the trace in~\eqref{eq:cubic-S0} is $3-33\omega$; the difference is $1+2\omega$.

\begin{proof}[Proof of Theorem~\ref{thm:B}(b)]
Theorem~\ref{thm:A}(c) gives $S_0(G_3)\subseteq3^9\Z$. For the reverse inclusion,
note that each row of the table has three character norms equal to $3$, one equal to
$9$, and $\Nm(\Delta(F_c))=3$, where $\Nm(r+s\omega)=r^2-rs+s^2$. Thus $Q(F_c)=3^8$
for $c=0,1,2,3,4$, and Lemma~\ref{lem:shift} gives
\begin{equation}\label{eq:G3-progressions}
 D_{G_3}(F_c+t\Sigma)=3^9(c+9t)\qquad(c=0,1,2,3,4,\ t\in\Z).
\end{equation}
Since $G_3$ has odd order, $D_{G_3}(-F)=-D_{G_3}(F)$. The residues
$0,\pm1,\pm2,\pm3,\pm4$ exhaust $\Z/9\Z$, so these progressions and their negatives
contain every multiple of $3^9$. The statement about valuations follows, because
$1\in S(G_3)$ and $3^e\in3^9\Z$ for every $e\ge9$.
\end{proof}

\subsection{Comparison with the other groups of order \texorpdfstring{$27$}{27}}

\begin{corollary}\label{cor:p3}\leavevmode
\begin{enumerate}
\item[(i)] A prime $q$ lies in $S(C_9\rtimes C_3)$ if and only if
$q\equiv\pm1\pmod{27}$.
\item[(ii)] $S^*(C_3\times C_9)\subsetneq S^*(C_9\rtimes C_3)$, and the smallest
positive integer in the difference is $62$.
\item[(iii)] The determinant sets of the five groups of order $27$ form the
strict inclusion chain
\[
 S(C_3^3)\subsetneq S(H_3)\subsetneq S(C_3\times C_9)
 \subsetneq S(C_9\rtimes C_3)\subsetneq S(C_{27}).
\]
\end{enumerate}
\end{corollary}

\begin{proof}
(i) If a prime $q$ lies in $S(G_3)$, then $q\in T_3$ by Corollary~\ref{cor:primes};
conversely, every $m\equiv\pm1\pmod{27}$ satisfies the criterion of
Theorem~\ref{thm:B}(a) with $n=1$.

(ii) By Theorem~3.4 of~\cite{MossinghoffPinner2026}, $S^*(C_3\times C_9)$ consists of
the integers $D\equiv\pm1\pmod{27}$ and the integers $D\equiv\pm8,\pm10\pmod{27}$ of
the form $D=A\Nm(B+C\omega)$ with $A\equiv\pm2,\pm4\pmod9$, $B\equiv A\pmod9$ and
$C\equiv A-B\pmod{27}$. The first kind lies in $S(G_3)$ by Theorem~\ref{thm:A}(b). For the second, write
$B=A+9s$. Then $C\equiv-9s\pmod{27}$ and
\[
 \Nm(B+C\omega)=B^2-BC+C^2\equiv A^2+18As+9As\equiv A^2\pmod{27},
\]
so $n=\Nm(B+C\omega)$ is a norm divisor of $D$ with
$Dn^3\equiv A\cdot A^8=A^9\equiv\pm1\pmod{27}$, because $(\Z/27\Z)^\times$ is cyclic
of order $18$; also $D\equiv A^3\equiv\pm1\pmod9$. So $D\in S(G_3)$ by
Theorem~\ref{thm:B}(a).

The inclusion is strict. On the one hand $62=2\cdot31$ lies in $S(G_3)$ with $n=31$,
since $62\cdot31^3\equiv-1\pmod{27}$. On the other hand $62\equiv8\pmod{27}$, so
$62\in S^*(C_3\times C_9)$ would require $\pm62=A\Nm(B+C\omega)$ with
$A\in\{\pm2,\pm31\}$, these being the divisors congruent to $\pm2$ or $\pm4$ modulo
$9$. Now $A=\pm31$ needs $\Nm(B+C\omega)=2$, which is not a norm, and $A=\pm2$ needs
$B^2-BC+C^2=31$ with $B\equiv A\pmod9$ and $C\equiv A-B\pmod{27}$; here
$31\ge\frac34\max(B^2,C^2)$ forces $|B|,|C|\le6$, hence $B=A$ and $C=0$, and then
$B^2-BC+C^2=4$. So $62\notin S^*(C_3\times C_9)$. Finally, by Theorem~\ref{thm:B}(a)
the positive elements of $S^*(G_3)$ below $62$ are $1$, $8$, $26$, $28$, $53$ and $55$
(the remaining elements of $T_2$ in this range have no suitable norm divisor; see the
list below). Of these, $8=2\Nm(2)$, and the others are congruent to $\pm1$ modulo
$27$, so all of them lie in $S^*(C_3\times C_9)$.

(iii) The first two inclusions are recorded in~\cite[\S1]{MossinghoffPinnerAbelian27}.
The coprime sets of $C_3^3$ and $H_3$ both equal $T_3$. Their divisible sets satisfy
$S_0(C_3^3)\subseteq3^{18}\Z$~\cite[Theorem~2.2]{MossinghoffPinner2026} and
$S_0(H_3)=3^{12}\Z$~\cite[Theorem~3.6]{MossinghoffPinnerHeisenberg}, so $3^{12}$
witnesses the first strict inclusion. The coprime set of $C_3\times C_9$ contains
$T_3$ and also $8=2\Nm(2)\notin T_3$ by the abelian criterion used in~(ii).
Its divisible set contains $3^{11}\Z$~\cite[Theorem~1(a)]{MossinghoffPinnerAbelian27},
so it contains $S_0(H_3)$. This gives the second strict inclusion.

For the third inclusion, the bound
$S_0(C_3\times C_9)\subseteq3^9\Z$~\cite[Theorem~3.2(iii)]{MossinghoffPinner2026}
and Theorem~\ref{thm:B}(b) give the inclusion on multiples of $3$. Part~(ii)
gives the strict inclusion on coprime values. The new abelian classification
also distinguishes the groups on multiples of $3$. At valuation $10$ it
requires a prime divisor of the cofactor~\cite[Theorem~1(b)]{MossinghoffPinnerAbelian27},
so $3^{10}$ lies in $S(G_3)\setminus S(C_3\times C_9)$.

Finally, $S(C_{27})$ contains all integers prime to $3$ and all multiples of $3^6$
\cite[Theorem~3.1]{MossinghoffPinner2023}. It therefore contains $S(G_3)$ by
Theorem~\ref{thm:B}(b). The integer $2$ belongs to $S(C_{27})$, and the congruence
in Theorem~\ref{thm:B}(a) excludes it from $S(G_3)$.
\end{proof}

For order $16$, by contrast, two nonisomorphic groups have the same set of integer
group determinants~\cite{BautistaPaudelPinner2024}.

Comparing the two criteria for all positive $m\le3000$, we find $71$ integers in
$S^*(G_3)\smallsetminus S^*(C_3\times C_9)$; those below $300$ are
$62,91,98,143,172,260,278,280,287$. The positive elements of $S(G_3)$ below $300$ that are prime to $3$ are
\begin{gather*}
 1,\ 8,\ 26,\ 28,\ 53,\ 55,\ 62,\ 64,\ 80,\ 82,\ 91,\ 98,\ 100,\ 107,\ 109,\ 116,\ 125,\ 134,\\
 136,\ 143,\ 161,\ 163,\ 172,\ 188,\ 190,\ 208,\ 215,\ 217,\ 224,\ 242,\ 244,\ 260,\ 269,\ 271,\\
 278,\ 280,\ 287,\ 296,\ 298 .
\end{gather*}
The elements of $T_2$ below $100$ that do not occur are
$10,17,19,35,37,44,46,71,73,89$. For instance, $91\cdot7^3\equiv1\pmod{27}$, so $91$
is admitted with $n=7$ (its other norm divisors $13$ and $91$ do not work), whereas
the norm divisors of $35=5\cdot7$ are $1$ and $7$, and neither $35$ nor $35\cdot7^3$
is congruent to $\pm1$ modulo $27$. The primes in $S(G_3)$ below $300$ are
$53,107,109,163,269,271$.

\section{Coprime determinants for \texorpdfstring{$p=5$}{p=5}}\label{sec:g5-coprime}

In this section and the next, $p=5$, $\cO=\Z[\omega]$ and $\pi=1-\omega$. We use the
isomorphism $\cO/5\cO\cong\F_5[w]/(w^4)$, $\omega\mapsto E(w)$, from the introduction
(it is one because $\Phi_5(z)\equiv(z-1)^4\pmod5$ and $E(w)-1$ is $w$ times a unit),
and the invariant $\ell_3$. The truncated logarithm $\log(1+z)=z-z^2/2+z^3/3$ modulo
$w^4$ has denominators prime to $5$ and turns products into sums, so $\ell_3$ is a
homomorphism from $(\cO/5\cO)^\times$ to $\F_5$. We call $\log\bigl(\alpha(w)/\alpha(0)\bigr)$ the normalized logarithm of
$\alpha$, and write $\ell_1(\alpha)$ and $\ell_2(\alpha)$ for its coefficients of $w$
and $w^2$.

By Theorem~\ref{thm:A}(b), the condition of Proposition~\ref{prop:nec} is
sufficient when it holds with $n=1$. For $p=5$ it is not sufficient in general, and
we first explain why.

\subsection{The cubic invariant}

For $f\in\Z[C_5^2]$ with augmentation $c=f(1,1)$ prime to $5$, put
\begin{equation}\label{eq:g5-binary-cubic}
 C_f(U,V)=\text{coefficient of $w^3$ in }
 \log\frac{f\bigl(E(Uw),E(Vw)\bigr)}{c}\ \in\F_5[U,V],
\end{equation}
where $c$ is reduced modulo $5$. The coefficient of $w^r$ in $E(Uw)$ and $E(Vw)$ is
homogeneous of degree $r$ in $U$ and $V$, so $C_f$ is a binary cubic form. Since
$E(iw)\equiv E(w)^i\pmod{w^4}$, its value at a nonzero pair $(i,j)$ is
\[
 C_f(i,j)=\ell_3\bigl(f(\omega^i,\omega^j)\bigr).
\]

\begin{lemma}\label{lem:cubic-obstruction}
Every global unit $\varepsilon\in\cO^\times$ has $\ell_3(\varepsilon)=0$, whereas
every $\alpha\in\cO$ with $\Nm(\alpha)=11$ has $\ell_3(\alpha)\ne0$. There is no
$f\in\Z[C_5^2]$ with augmentation prime to $5$ such that, of the six norms
$\Nm(f(\omega^i,\omega^j))$, one for each line through the origin in $\F_5^2$, five
equal $1$ and one equals $11$.
\end{lemma}

\begin{proof}
Let $\varepsilon$ be a unit. The quotient $\varepsilon/\bar\varepsilon$ is a root of
unity~\cite[Ch.~1]{Washington1997}, so it equals $\pm\omega^r$ for some $r$, and
reduction modulo $\pi$ shows that it is $\omega^r$. Complex conjugation acts on
$\cO/5\cO$ by $w\mapsto-w$, so the normalized logarithms of $\varepsilon(w)$ and
$\varepsilon(-w)$ differ by $rw$. Comparing coefficients of $w^3$ gives
$2\ell_3(\varepsilon)=0$.

The four elements $2+\omega^r$, $1\le r\le4$, have norm $11$. They generate distinct
prime ideals, because reduction modulo $(2+\omega^r)$ gives $\omega^r=-2$ in
$\F_{11}$, and the four resulting values of $\omega$ are distinct primitive fifth
roots of unity. So these are all the primes above $11$, and every element of norm
$11$ is a unit times one of the four. Now
\[
 \ell_3(2+\omega^r)=\text{coefficient of $w^3$ in }\log\frac{2+E(rw)}{3}
 =\frac{r^3}{81}=r^3\ne0\quad\text{in }\F_5 ,
\]
and the second assertion follows from additivity and the case of units.

For the last assertion, suppose that five of these norms equal $1$. The
corresponding values are units, so $C_f$ vanishes at five points of $\mathbf P^1(\F_5)$. A
nonzero binary cubic form has at most three zeros on the projective line, so $C_f=0$.
This contradicts the fact that the value of norm $11$ on the remaining orbit has
nonzero $\ell_3$.
\end{proof}

\begin{proposition}\label{prop:418}
The integer $418$ satisfies the condition of Proposition~\ref{prop:nec} for $p=5$,
with norm divisor $n=11$, but $418\notin S(G_5)$.
\end{proposition}

\begin{proof}
We have $11=\Nm(2+\omega)$, and
\[
 418^4\equiv1\pmod{25},\qquad
 (418\cdot11^5)^4\equiv1\pmod{125},\qquad
 418^4\equiv51\pmod{125},
\]
so $418$ meets the condition and lies outside $T_3$. Suppose that
$D_{G_5}(F)=418$. Let $c=F(1)$, and choose one value $\alpha_L$ in each of the six
nontrivial orbits of linear characters, indexed by the points $L$ of
$\mathbf P^1(\F_5)$. The Frobenius factorization gives
\[
 418=c\prod_{L\in\mathbf P^1(\F_5)}\Nm(\alpha_L)\,\Nm(\Delta(F))^5 ,
\]
with all norm factors positive integers. Since $418=2\cdot11\cdot19$ is squarefree,
$\Nm(\Delta(F))=1$. The residue degrees of $2$ and $19$ in $\Q(\omega)$ are $4$ and
$2$, so neither can occur to the first power in a character norm. Hence either all six
character norms are $1$, or five are $1$ and one is $11$. The second case is excluded
by Lemma~\ref{lem:cubic-obstruction}. In the first case $c=418$ and $A(F)=418$, so
Theorem~\ref{thm:A}(a) would give $418\equiv418^5\pmod{125}$; but the two sides are
$43$ and $68$ modulo $125$.
\end{proof}

\subsection{The coprime criterion}\label{sec:g5-complete}

We now show that the cubic invariant is the only further obstruction. The
reconstruction of linear character values with prescribed norms and
$\ell_3$-invariants is modelled on~\cite[Lemmas~2.2, 3.1 and~3.2]{Panraksa2026Bernoulli},
where similar coordinates are used for $C_p\times C_p$, and prime-factor tests of the
kind given in Section~\ref{sec:g5-prime-test} also appear
in~\cite[Theorem~2.3]{Panraksa2026Order24}. What is new here is the compatibility of
these values with $\Delta$.

\begin{lemma}\label{lem:g5-units}
The image of $\cO^\times$ in $(\cO/\pi^5\cO)^\times$ is
\[
 \{a:\Nm(a)\equiv1\pmod{25},\ \ell_3(a)=0\},
\]
and it consists of the $100$ distinct classes $\omega^r(1+\omega)^s$, $0\le r<5$,
$0\le s<20$. The map $\varepsilon\mapsto(\varepsilon\bmod\pi,\ \ell_1(\varepsilon),\
\ell_2(\varepsilon))$ sends $\cO^\times$ onto $\F_5^\times\times\F_5^2$. In particular, if $\ell_3(\alpha)=0$ and $5\nmid c$, then some global
unit $\varepsilon$ satisfies $\varepsilon\alpha\equiv c\pmod{5\cO}$.
\end{lemma}

\begin{proof}
There are $2500$ units modulo $\pi^5$. The norm modulo $25$ is well defined on these
classes: if $z\in\pi^5\cO_\p=5\pi\cO_\p$, then $\Tr(z)\in25\Z_5$ and all higher terms
of $\Nm(1+z)$ are divisible by $25$. So
\[
 a\longmapsto\bigl(\Nm(a)\bmod25,\ \ell_3(a)\bigr)
\]
is a homomorphism to $(1+5\Z/25\Z)\times\F_5$. It is surjective, since $6$ maps to
$(21,0)$ and $1+\pi^3$ has nonzero $\ell_3$, and so its kernel has order $100$. Every
global unit lies in the kernel: its norm is $1$, since $k$ is totally complex, and its
invariant $\ell_3$ vanishes by Lemma~\ref{lem:cubic-obstruction}.

To find $100$ distinct classes of global units, we use $\ell_1$ and $\ell_2$. The
identity
\[
 \log\frac{1+E(w)}{2}\equiv\frac w2+\frac{w^2}{8}\pmod{w^4}
\]
shows that $\ell_2(1+\omega)=2$, while $\ell_2(\omega)=0$. If two of the classes
$\omega^r(1+\omega)^s$ coincide, reduction modulo $\pi$ shows that $s-s'$ is divisible
by $4$, and applying $\ell_2$ shows that it is divisible by $5$. Hence $s=s'$, and
then $\ell_1(\omega)=1$ gives $r=r'$. These classes therefore fill the kernel. The same
argument shows that their triples (residue modulo $\pi$, $\ell_1$, $\ell_2$) run
through all of $\F_5^\times\times\F_5^2$. Now choose $\varepsilon$ so that
$\varepsilon\alpha$ has triple $(c,0,0)$; since also $\ell_3(\varepsilon\alpha)=0$, the
image of $\varepsilon\alpha$ in $\F_5[w]/(w^4)$ is the constant $c$.
\end{proof}

The next lemma describes $\Delta$ modulo $5$. Part~(b) is needed only in
Section~\ref{sec:g5-divisible}.

\begin{lemma}\label{lem:g5-delta}
\leavevmode
\begin{enumerate}
\item[(a)] For every $F\in\Z_5[G_5]$,
\[
 \Delta(F)\equiv\Fbar(\omega,1)\pmod{5\cO_\p}.
\]
\item[(b)] For all $A,B\in\Lam_\p$, every coefficient of the polynomial
$\det(A+tB)-\det A-t^5\det B\in\cO_\p[t]$ lies in $5\cO_\p$. In particular,
\[
 \det(A+\pi B)\equiv\det A\pmod{\pi^5\cO_\p}.
\]
\end{enumerate}
\end{lemma}

\begin{proof}
(a) Suppose first that $F\in\Z[G_5]$ has augmentation prime to $5$. Both sides of the
congruence define multiplicative maps from $R^\times=\F_5[C_5^2]^\times$ to
$(\cO/5\cO)^\times$. For the left-hand side, two lifts of the same element of $R$
differ by an element of $5\Z_5[G_5]+(x^5-1)\Z_5[G_5]$. A lift with unit augmentation
is a unit of $\Z_5[G_5]$, so their images under $\rho$ differ multiplicatively by an
element of $1+\pi\Lam_\p$, and by Lemma~\ref{lem:newton} their determinants agree
modulo $\pi^5$, hence modulo $5$. It remains to compare the two maps on the generators
of Lemma~\ref{lem:basis}. For a constant $c\in\F_5^\times$ both give $c^5=c$. For
$e_{ij}=1+(1-X)^i(1-Y)^j$, put $B=(1-\rho(x))^i(1-u)^j$. The Newton-identity argument
of Lemma~\ref{lem:newton} gives $\det(1+B)\equiv1+\det B\pmod{5\cO}$. If $j>0$ then
$\det B=0$, and if $j=0$ then $i>0$ and $\det B=\Nm_{K/k}(1-\zeta)^i=\pi^i$. These are
exactly the values of $e_{ij}(\omega,1)$.

Now let $F\in\Z[G_5]$ be arbitrary, and choose $h\in\Z$ such that $F+h$ has
augmentation prime to $5$. Lemma~\ref{lem:newton} gives
$\det(\rho(F)+hI)\equiv\Delta(F)+h^5\pmod{5\cO}$, and $h^5\equiv h\pmod5$, so the case
already proved, applied to $F+h$, gives $\Delta(F)\equiv\Fbar(\omega,1)\pmod{5\cO}$.
Both sides are continuous in the coefficients of $F$, so the congruence extends to
$\Z_5[G_5]$.

(b) The image of $\Z_5[G_5]$ under $\rho$ is $\Lam_\p$, since the central elements
$x^{5j}$ supply the coefficients in $\cO_\p$. As $F\mapsto\Fbar(\omega,1)$ is additive,
(a) gives
\[
 \det(A+hB)\equiv\det A+h\det B\pmod{5\cO_\p}\qquad(h\in\Z).
\]
Write $\det(A+tB)=\det A+\sum_{j=1}^4c_jt^j+t^5\det B$. Since $h^5\equiv h\pmod5$,
the polynomial $\sum_{j=1}^4c_jt^j$ vanishes modulo $5\cO_\p$ at $t=1,2,3,4$, and as
$\det(t^j)_{1\le t,j\le4}=288$ is prime to $5$, we get $c_j\in5\cO_\p$ for
$1\le j\le4$.
Substituting $t=\pi$ gives the last assertion, because $5\pi\cO_\p=\pi^5\cO_\p$.
\end{proof}

\begin{lemma}\label{lem:g5-reconstruct}
Let $c\in\Z$ be prime to $5$, and let $n_L\in\cN$ be six integers indexed by
$L\in\cL=\{(b,1):b\in\F_5\}\cup\{(1,0)\}$. If $c\prod_Ln_L\in T_2$, then there is an
$f\in\Z[C_5^2]$ with augmentation $c$ such that
\[
 \Nm\bigl(f(\omega^{L_1},\omega^{L_2})\bigr)=n_L,
 \qquad \ell_3\bigl(f(\omega^{L_1},\omega^{L_2})\bigr)=0
 \qquad(L\in\cL).
\]
\end{lemma}

\begin{proof}
Choose $\alpha_L\in\cO$ with norm $n_L$ and $\ell_3(\alpha_L)=0$. By
Lemma~\ref{lem:g5-units}, after multiplying by global units we may assume that every
$\alpha_L$ is congruent to $c$ modulo $5\cO$; this changes neither the norm nor
$\ell_3$. Write $\alpha_L=c+5\beta_L$. Expanding the norms gives
\[
 n_L\equiv c^4+5c^3\Tr(\beta_L)\pmod{25},\qquad
 c\prod_Ln_L\equiv c^{25}\Bigl(1+\frac5c\sum_L\Tr(\beta_L)\Bigr)\pmod{25}.
\]
Since $c^{100}\equiv1\pmod{25}$, the hypothesis implies
$\sum_L\Tr(\beta_L)\equiv0\pmod5$.

Prescribe the trivial character value to be $c$, the value at $L$ to be $\alpha_L$,
and the value at $aL$ to be $\tau_a(\alpha_L)$. Fourier inversion gives the rational
coefficients
\[
 b_{ij}=\frac1{25}\Bigl(c+\sum_{L\in\cL}
 \Tr\bigl(\alpha_L\omega^{-L_1i-L_2j}\bigr)\Bigr),\qquad 0\le i,j<5 .
\]
The terms involving $c$ (the first term and the summand $c$ of each
$\alpha_L=c+5\beta_L$)
contribute $25c$ to the numerator for $(i,j)=(0,0)$ and nothing otherwise, because a nonzero vector is orthogonal to exactly
one line and $\Tr(\omega^t)$ is $4$ or $-1$ according as $t=0$ or $t\ne0$. The rest of
the numerator is $5\sum_L\Tr(\beta_L\omega^{-L_1i-L_2j})$. For every $\beta\in\cO$ and
every integer $t$ we have $\Tr(\beta\omega^t)\equiv\Tr(\beta)\pmod5$, since the
difference is the trace of an element of $\pi\cO$. Hence this part of the numerator is
also divisible by $25$. All $b_{ij}$ are integers, and $f=\sum b_{ij}X^iY^j$ has the
prescribed values.
\end{proof}

\begin{proof}[Proof of Theorem~\ref{thm:C}(a)]
Suppose first that $m=D_{G_5}(F)$ is prime to $5$. Put $f=\Fbar$, $c=f(1,1)$,
$\alpha_b=f(\omega^b,\omega)$ for $b\in\F_5$ and $\alpha_\infty=f(\omega,1)$, and set
\[
 n=\prod_{b\in\F_5}\Nm(\alpha_b),\qquad
 h=\Nm(\alpha_\infty),\qquad s=\Nm(\Delta(F)).
\]
Then $m=chns^5$ and $A(F)=ch$; in particular $n$ is a positive divisor of $|m|$. The
binary cubic form $C_f$ of~\eqref{eq:g5-binary-cubic} gives
\[
 \ell_3\Bigl(\prod_b\alpha_b\Bigr)=\sum_{b\in\F_5}C_f(b,1)=0 ,
\]
because $\sum_bb^j=0$ in $\F_5$ for $0\le j\le3$. Hence $n\in\cN$.
Theorem~\ref{thm:A}(a) gives $m\equiv(ch)^5\pmod{125}$. Every norm prime to $5$ is
$1$ modulo $5$, so $s^{25}\equiv1\pmod{125}$, and substituting $ch=m/(ns^5)$ gives
\begin{equation}\label{eq:g5-scalar}
 m^4\equiv n^5\pmod{125}.
\end{equation}
Since $n^{25}\equiv1\pmod{125}$, this is equivalent to $mn^5\in T_3$. Also $m\in T_2$
by Proposition~\ref{prop:nec}.

Conversely, let $m\in T_2$, and let $n\in\cN$ divide $|m|$ with $mn^5\in T_3$. Put
$c=m/n$. Lemma~\ref{lem:g5-reconstruct}, with $n_{(0,1)}=n$ and the other five norms
equal to $1$, gives $f$ with $M_1(f)=m$, $A(f)=c$, and all six values having
$\ell_3=0$. Choose any lift $F'$ of $f$ and put $\delta_0=\Delta(F')$. By
Lemma~\ref{lem:g5-delta}(a), $\ell_3(\delta_0)=\ell_3(f(\omega,1))=0$.
Congruence~\eqref{eq:g5-scalar} is equivalent to $m\equiv c^5\pmod{125}$, so
Theorem~\ref{thm:A}(a) and the Frobenius factorization give
$\Nm(\delta_0)^5\equiv1\pmod{125}$; as $\Nm(\delta_0)\equiv1\pmod5$, this means
$\Nm(\delta_0)\equiv1\pmod{25}$. Lemma~\ref{lem:g5-units} therefore supplies a global
unit $\varepsilon$ with $\varepsilon\equiv\delta_0\pmod{\pi^5}$. By
Lemma~\ref{lem:L1} and Proposition~\ref{prop:fibre}, there is a lift $F$ of $f$ with
$\Delta(F)=\varepsilon$, and then $D_{G_5}(F)=M_1(f)\Nm(\varepsilon)^5=m$.
\end{proof}

\subsection{A prime-factor test}\label{sec:g5-prime-test}

Membership in $\cN$ can be decided from the prime factorization. The ring $\cO$ has
class number one~\cite[Theorem~11.1]{Washington1997}, so for an ideal $I$ of $\cO$
prime to $5$ we may put $d(I)=\ell_3(\alpha)$, where $I=(\alpha)$. By
Lemma~\ref{lem:cubic-obstruction} this does not depend on the generator. It is
additive on products, and $d(\tau_aI)=a^3d(I)$. (These are the cubic ideal
coordinates of~\cite[Lemma~3.2]{Panraksa2026Bernoulli}.) In this notation, a prime
$q\equiv1\pmod5$ lies in $\cB$ if and only if $d(\mathfrak q)\ne0$ for a prime ideal
$\mathfrak q$ above $q$; since the prime ideals above $q$ are the $\tau_a\mathfrak q$,
the choice of $\mathfrak q$ does not matter.

\begin{proposition}\label{prop:g5-norms}
Let $n=\prod_{q\ne5}q^{e_q}$ be a positive integer, and put $f_q=\ord_5(q)$, the
order of $q$ modulo $5$. Then $n\in\cN$ if and only if
\[
 f_q\mid e_q\ \text{ for every }q,
 \qquad\text{and}\qquad \sum_{q\in\cB}e_q\ne1 .
\]
\end{proposition}

\begin{proof}
A prime ideal $\mathfrak q$ above $q$ has norm $q^{f_q}$, and its decomposition group
contains $\tau_q$, so $(q^3-1)d(\mathfrak q)=0$. If $q\not\equiv1\pmod5$, then
$q^3\ne1$ in $\F_5^\times$, and so $d(\mathfrak q)=0$. If $q\equiv1\pmod5$, the four
values $d(\tau_a\mathfrak q)=a^3d(\mathfrak q)$ are either all zero or, since cubing
permutes $\F_5^\times$, run through $\F_5^\times$. So each prime factor of $n$ in
$\cB$, counted with multiplicity, contributes a nonzero element of $\F_5$ that we may
choose freely, and all other prime-ideal factors contribute zero. A single nonzero
element of $\F_5$ is not zero, but for any number $N\ge2$ of them there is a choice
with sum zero: take pairs $1,4$, together with $1,1,3$ if $N$ is odd. This proves the
assertion for norms of ideals, and since $\cO$ is principal, these ideals give
elements of the same norm.
\end{proof}

Primes $q\equiv1\pmod5$ can be tested for membership in $\cB$ without finding a
generator of a prime ideal above $q$. Choose a primitive root $g$ modulo $q$, write
$\operatorname{ind}_g(z)$ for the discrete logarithm modulo $5$, and put
\[
 A_j=\#\{z\in\F_q\smallsetminus\{0,1\}:\operatorname{ind}_g\bigl(z(1-z)\bigr)=j\},
 \qquad S_i=\sum_{j=0}^4j^iA_j\bmod5\quad(1\le i\le3).
\]

\begin{proposition}\label{prop:g5-jacobi}
A prime $q\equiv1\pmod5$ lies outside $\cB$ if and only if
\[
 S_3+3S_1S_2+2S_1^3=0\quad\text{in }\F_5 .
\]
In particular, the test does not depend on the choice of $g$.
\end{proposition}

\begin{proof}
Let $\chi$ be the character of $\F_q^\times$ with $\chi(g)=\omega$, extended by
$\chi(0)=0$, and let $J=\sum_z\chi(z)\chi(1-z)=\sum_jA_j\omega^j$ be its Jacobi sum.
The conjugates of $J$ are the Jacobi sums $J(\chi^a,\chi^a)$, $1\le a\le4$, and
since $\chi^a$ and $\chi^{2a}$ are nontrivial, each has absolute value
$\sqrt q$~\cite[Ch.~8, \S3]{IrelandRosen1990}; so $\Nm(J)=q^2$. Put $r=(q-1)/5$, and
let $\mathfrak q_b$ be the prime ideal above $q$ defined by $\omega\mapsto g^{br}$.
Modulo $\mathfrak q_b$, $J$ reduces to $\sum_zz^{br}(1-z)^{br}$. For $b=1,2$ this
vanishes, since the degree $2br$ is less than $q-1$. For $b=3,4$ the only positive
multiple of $q-1$ among the exponents is $q-1$ itself, and its binomial coefficient is
nonzero modulo $q$. So $J$ lies in $\mathfrak q_1$ and $\mathfrak q_2$ and in no other
prime above $q$, and comparing norms gives $(J)=\mathfrak q_1\mathfrak q_2$. As
$\mathfrak q_2=\tau_3\mathfrak q_1$, we get
$\ell_3(J)=(1+3^3)d(\mathfrak q_1)=3d(\mathfrak q_1)$ in $\F_5$.

Finally, $\sum_jA_j=q-2\equiv-1\pmod5$, and expanding the normalized logarithm gives
\[
 \ell_3(J)=-S_3+2S_1S_2-2S_1^3=-(S_3+3S_1S_2+2S_1^3)\quad\text{in }\F_5 . \qedhere
\]
\end{proof}

Among the $22$ primes $q\equiv1\pmod5$ below $500$, only $211$, $281$, $421$ and
$461$ lie outside $\cB$. By Proposition~\ref{prop:g5-norms}, $11\notin\cN$, while $11^2$, $11^3$ and the product of any
two primes in $\cB$ lie in $\cN$. For example, Theorem~\ref{thm:C}(a) admits $2526$
with $n=421$, $3782$ with $n=31\cdot61$, and $3993$ with $n=11^3$. The positive
elements of $S^*(G_5)$ below $250$ are
\[
 1,\ 32,\ 57,\ 68,\ 124,\ 126,\ 176,\ 182,\ 193,\ 243,\ 249,
\]
and there are $365$ of them up to $10000$. The norm-divisor condition of
Proposition~\ref{prop:nec} admits $84$ further integers in that range, the first five
being $418,451,524,549,574$.

\section{Divisible determinants for \texorpdfstring{$p=5$}{p=5}}\label{sec:g5-divisible}

We keep the notation of Section~\ref{sec:g5-coprime}. The proof of
Theorem~\ref{thm:C}(b) needs $\Delta(F)$ modulo $\pi^5$, one power of $\pi$ more than
Lemma~\ref{lem:g5-delta}(a) gives. We first show that if the image of $F$ in
$\F_5[C_5\times C_5]$ has a nonzero linear part, then $\Delta(F')$ runs through all of
$\Delta(F)+\pi^5\cO_\p$ as $F'$ runs over the lifts of $\Fbar$ to $\Z_5[G_5]$. We then
express $\ell_3(\Delta(F)/\pi)$ in terms of the linear character values.

\subsection{Values of \texorpdfstring{$\Delta$}{Delta} modulo \texorpdfstring{$\pi^5$}{pi\^{}5}}

\begin{lemma}\label{lem:g5-singular}
Identify $\Lam_\p/\pi\Lam_\p$ with the ring $R=\F_5[e,v]/(e^5,v^5)$ of
Section~\ref{sec:order} via $e=1-x$ and $v=1-y$, and let $\m=(e,v)$. If $A\in\Lam_\p$
has image $a\in\m\smallsetminus\m^2$, then
\[
 \nrd(A+\pi\Lam_\p)=\nrd(A)+\pi^5\cO_\p .
\]
The hypothesis on $a$ holds whenever $v_\pi(\nrd A)=1$.
\end{lemma}

\begin{proof}
We use a second model of $\Lam_\p$. Let $\mathsf T$ be multiplication by $\zeta$ on
the ring of integers of $\Q_5(\zeta)$ and $\mathsf U$ the action of $\sigma$, both
written in the $\cO_\p$-basis $1,\zeta,\dots,\zeta^4$. Then
\[
 \mathsf U=\diag(1,\omega,\omega^2,\omega^3,\omega^4),\qquad
 \mathsf U\mathsf T=\omega\mathsf T\mathsf U,\qquad
 \mathsf T^5=\omega I,\qquad \mathsf U^5=I,
\]
and $x\mapsto\mathsf T$, $y\mapsto\mathsf U$ maps $\Z_5[G_5]$ onto
$\Lam_\p=\bigoplus_{0\le i,j<5}\cO_\p\mathsf T^i\mathsf U^j$. This representation is
irreducible of degree $5$ with $x^5\mapsto\omega$, hence equivalent to $\rho$, so the
determinant is again the reduced norm; and the relations show that
$\Lam_\p/\pi\Lam_\p=R$.

The inclusion $\nrd(A+\pi\Lam_\p)\subseteq\nrd(A)+\pi^5\cO_\p$ is
Lemma~\ref{lem:g5-delta}(b). For the reverse inclusion, put
$\trz(C)=\tr(C)/5$ for $C\in\Lam_\p$. The trace of $\mathsf T^i\mathsf U^j$ is zero
unless $i=j=0$, when it is $5$, so the pairing $(C,B)\mapsto\trz(CB)$ is perfect:
each basis element pairs nontrivially only with its inverse, up to a unit of $\cO_\p$.
Its reduction is a perfect associative pairing on $R$. The annihilator of $\m$ in $R$
is the socle $\F_5e^4v^4$, and hence
\[
 \trz(c\m)=0\iff c\in\F_5e^4v^4 .
\]
Indeed, if $\trz(c\m)=0$, then for $h\in\m$ and $r\in R$ we have $\trz(chr)=0$, since
$hr\in\m$, so $ch=0$ by perfectness; the converse is clear.

Lemma~\ref{lem:newton} and the Cayley--Hamilton theorem give
$\adj(A)\equiv A^4\pmod{5\Lam_\p}$. Since $a$ has a nonzero linear part, $a^4$ has a
nonzero homogeneous term of degree $4$, so it is not in the socle, which has degree
$8$. We may therefore choose $B\in\Lam_\p$ with image in $\m$ such that
$\trz(\adj(A)B)$ is a unit. The image of $B$ has zero augmentation, so
$\det B\in\pi\cO_\p$ by Lemma~\ref{lem:g5-delta}(a). Consider
$\phi(t)=\det(A+\pi tB)$. Its coefficient of $t$ is
$\pi\tr(\adj(A)B)=5\pi\trz(\adj(A)B)$, of valuation exactly $5$. By
Lemma~\ref{lem:g5-delta}(b) the coefficients of $t^2,t^3,t^4$ have valuations at
least $6,7,8$, and the coefficient of $t^5$ is $\pi^5\det B$, of valuation at least
$6$. Hence
\[
 \phi(t)=\det A+\pi^5\bigl(\varepsilon_0t+\pi\phi_1(t)\bigr),
 \qquad \varepsilon_0\in\cO_\p^\times,\quad \phi_1\in\cO_\p[t],\quad \phi_1(0)=0 .
\]
For each $z\in\cO_\p$ the equation $\varepsilon_0t+\pi\phi_1(t)=z$ has a simple root modulo
$\pi$, so Hensel's lemma gives a solution $t\in\cO_\p$. This proves the reverse
inclusion.

Finally, suppose that $v_\pi(\det A)=1$, and let $F\in\Z_5[G_5]$ map to $A$. By
Lemma~\ref{lem:g5-delta}(a), $\Fbar(\omega,1)$ has valuation $1$. Write the
reduction of $\Fbar$ modulo $5$ as a polynomial in $e$ and $v$ and substitute
$e=\pi$, $v=0$: its constant term must vanish and its coefficient of $e$ must not.
Hence $a\in\m\smallsetminus\m^2$.
\end{proof}

\subsection{Two congruences}

For $f\in\Z[C_5^2]$ let $c=f(1,1)$, and let $\theta(f)$ be the coefficient of the
identity in $f^5$. We index the six nontrivial character orbits by
\[
 L_\infty=(1,0),\qquad L_a=(a,1)\quad(0\le a\le4),
\]
and write $\alpha_L=f(\omega^{L_1},\omega^{L_2})$; thus
$\alpha_\infty=f(\omega,1)$.

\begin{lemma}\label{lem:g5-divisible-scalar}
For every $F\in\Z[G_5]$ with $\Fbar=f$,
\begin{equation}\label{eq:g5-divisible-scalar}
 \Delta(F)\equiv\alpha_\infty+\theta(f)-c\pmod{\pi^5\cO}.
\end{equation}
Moreover $5\mid\theta(f)-c$. If $v_\pi(\alpha_\infty)=1$ and $a_\infty$ is the residue
of $\alpha_\infty/\pi$ modulo $\pi$, then
\begin{equation}\label{eq:g5-normalized-cubic}
 \ell_3\bigl(\Delta(F)/\pi\bigr)=\ell_3(\alpha_\infty/\pi)
 +\frac{(\theta(f)-c)/5}{a_\infty}\quad\text{in }\F_5 .
\end{equation}
\end{lemma}

\begin{proof}
By Lemmas~\ref{lem:lifts} and~\ref{lem:g5-delta}(b), $\Delta(F)$ modulo $\pi^5$
depends only on $f$. We may therefore take $F=\sum b_{ij}x^iy^j$, where
$f=\sum_{i,j=0}^4b_{ij}X^iY^j$; in the model of Lemma~\ref{lem:g5-singular} its image
is $A=\sum b_{ij}\mathsf T^i\mathsf U^j$. Write $s_j=\tr(A^j)=5t_j$. Newton's identity
for the fifth elementary symmetric function reads
\[
 120\det A=s_1^5-10s_1^3s_2+15s_1s_2^2+20s_1^2s_3-20s_2s_3-30s_1s_4+24s_5 .
\]
After substituting $s_j=5t_j$, every term except $24s_5$, divided by $120$, lies in
$25\cO_\p$. Hence $\det A\equiv\trz(A^5)\pmod{25\cO_\p}$.

Expand $A^5$ into words of length five in the matrices $\mathsf T^i\mathsf U^j$.
Cyclic rotation divides the nonconstant words into orbits of size five, with the same
coefficient product and the same trace throughout an orbit. A word has nonzero trace
exactly when its total exponents of $\mathsf T$ and of $\mathsf U$ are divisible by
five, and then $\trz$ of the word is a power of $\omega$. So, modulo $5\pi\cO$, each
such orbit contributes its integer coefficient times $5$, and these are precisely the
nonconstant contributions to $\theta(f)$. A constant word contributes
$b_{ij}^5\omega^i$, because $(\mathsf T^i\mathsf U^j)^5=\omega^iI$, whereas its
contribution to $\theta(f)$ is $b_{ij}^5$. It follows that
\[
 \trz(A^5)\equiv\theta(f)+\sum_{i,j}b_{ij}^5(\omega^i-1)
 \equiv\theta(f)+\alpha_\infty-c\pmod{5\pi\cO},
\]
the second congruence because $5\mid b_{ij}^5-b_{ij}$. Since
$25\cO\subseteq5\pi\cO=\pi^5\cO$, this proves~\eqref{eq:g5-divisible-scalar}. Also
$f^5\equiv c\pmod{5\Z[C_5^2]}$, so $\theta(f)\equiv c\pmod5$.

For~\eqref{eq:g5-normalized-cubic}, note that $\pi^4\equiv-5\pmod{\pi^5\cO}$, so
$5/\pi\equiv-\pi^3\pmod{\pi^4\cO}$, and $-\pi^3\equiv w^3\pmod{w^4}$ under
$\omega=E(w)$. Dividing~\eqref{eq:g5-divisible-scalar} by $\pi$ therefore adds
$\bigl((\theta(f)-c)/5\bigr)w^3$ to the unit $\alpha_\infty/\pi$, and this changes the
coefficient of $w^3$ in its normalized logarithm by $\bigl((\theta(f)-c)/5\bigr)/a_\infty$.
\end{proof}

\begin{lemma}\label{lem:g5-fifth-trace}
For every $\gamma\in\cO_\p^\times$ with residue $a\in\F_5^\times$,
\begin{equation}\label{eq:g5-fifth-trace}
 \frac{\Tr\bigl((\pi\gamma)^5\bigr)}{125}\equiv a\,\ell_3(\gamma)\pmod5,
 \qquad
 \Tr\bigl((\pi^2\gamma)^5\bigr)\equiv0\pmod{625}.
\end{equation}
Here the first quotient lies in $\Z_5$, and $\Tr$ and $\ell_3$ are extended to the
completion.
\end{lemma}

\begin{proof}
Write $\gamma\equiv a+b\pi+c\pi^2+d\pi^3\pmod{\pi^4\cO_\p}$ with $a\in\{1,2,3,4\}$
and $b,c,d\in\{0,1,2,3,4\}$. Changing $\gamma$ by $\pi^4\beta$ changes
$(\pi\gamma)^5$ by an element of $\pi^{13}\cO_\p$: the first binomial term has
valuation at least $13$, and the others are no smaller. Since
$\pi^{13}\cO_\p=125\pi\cO_\p$ and $\Tr(\pi\cO_\p)\subseteq5\Z_5$, this changes the
trace by a multiple of $625$.

By~\eqref{eq:c}, $\Tr(\pi^n)=5\kappa(n)$ for $n\ge1$, and for $n=5,6,7,8,9,10$
these traces are $0$, $-25$, $-100$, $-275$, $-625$ and $-1250$. In the expansion of
$\pi^5(a+b\pi+c\pi^2+d\pi^3)^5$, the term $a^5\pi^5$ has trace $0$, and the other terms
that can contribute to the trace modulo $625$ are
\begin{align*}
 &5a^4\bigl(b\Tr(\pi^6)+c\Tr(\pi^7)+d\Tr(\pi^8)\bigr)\\
 &\quad+10a^3\bigl(b^2\Tr(\pi^7)+2bc\Tr(\pi^8)\bigr)+10a^2b^3\Tr(\pi^8).
\end{align*}
Every omitted mixed term has a coefficient divisible by $5$ and $\pi$-degree at least
$9$, and $\Tr(\pi^9\cO_\p)\subseteq125\Z_5$. The pure $b^5$ term has trace
$b^5\Tr(\pi^{10})$, divisible by $625$, and the pure $c^5$ and $d^5$ terms lie in
$\pi^{15}\cO_\p$, whose trace lies in $625\Z_5$. So the trace is divisible by $125$,
and the quotient reduces to
\[
 -b+c-d+\frac{2b^2+bc}{a}+\frac{3b^3}{a^2}.
\]
On the other hand, substituting $\pi=-w-w^2/2-w^3/6\pmod{w^4}$ in the normalized
logarithm gives
\[
 a\,\ell_3(\gamma)=-\frac b6+c-d-\frac{b^2}{2a}+\frac{bc}{a}-\frac{b^3}{3a^2},
\]
and the two expressions agree in $\F_5$. This proves the first congruence.

For the second, write $\gamma=a+z$ with $z\in\pi\cO_\p$. The constant term in
$\pi^{10}(a+z)^5$ has trace $a^5\Tr(\pi^{10})\in625\Z_5$. Each intermediate binomial
term lies in $5\pi^{11}\cO_\p$, whose trace lies in $625\Z_5$ because
$\pi^{11}\cO_\p=25\pi^3\cO_\p$. The last term lies in $\pi^{15}\cO_\p$ and has the same
property.
\end{proof}

\subsection{Proof of Theorem~\ref{thm:C}(b)}

We first show that a group determinant of valuation $13$ has a prime factor in $\cB$.
Suppose that $D_{G_5}(F)=5^{13}m$ with $5\nmid m$, and put $f=\Fbar$ and $c=f(1,1)$.
The proof of Theorem~\ref{thm:A}(c) gives $v_5(M_1(F))\ge8$ and
$v_\pi(\Delta(F))\ge1$, and since $D_{G_5}(F)=M_1(F)\Nm(\Delta(F))^5$, both are
equalities. Write the reduction of $f$ modulo $5$ as a polynomial in $e=1-X$ and
$v=1-Y$. Its linear part is nonzero, for otherwise each of the six values $\alpha_L$
would have valuation at least $2$, and $v_5(M_1(F))\ge1+12=13$. A nonzero linear form
vanishes on exactly one line of $\F_5^2$, so five of the $\alpha_L$ have valuation
exactly $1$, and the remaining one, say $\alpha_{L_*}$, has valuation at least $2$.
Since $v_5(c)\ge1$ and $v_5(M_1(F))=8$, the value $\alpha_{L_*}$ has valuation exactly
$2$ and $v_5(c)=1$. By Lemma~\ref{lem:g5-delta}(a), $\alpha_\infty\equiv\Delta(F)
\pmod5$ has valuation $1$, so $L_*\ne L_\infty$.

Write $c=5r$, $\alpha_L=\pi\eta_L$ for $L\ne L_*$, $\alpha_{L_*}=\pi^2\eta_{L_*}$ and
$\Delta(F)=\pi\eta$, where $\eta_L$, $\eta_{L_*}$ and $\eta$ lie in $\cO$ and are prime
to $\pi$. For $L\ne L_*$ put $a_L=\eta_L\bmod\pi$ and $d_L=\ell_3(\eta_L)$. By the
orthogonality of characters,
\[
 25\,\theta(f)=c^5+\sum_{L\in\mathbf P^1(\F_5)}\Tr(\alpha_L^5),
\]
and Lemma~\ref{lem:g5-fifth-trace} gives
\begin{equation}\label{eq:g5-boundary-trace}
 \frac{\theta(f)}5\equiv\sum_{L\ne L_*}a_Ld_L\pmod5 .
\end{equation}
Together with~\eqref{eq:g5-normalized-cubic}, this yields
\begin{equation}\label{eq:g5-boundary-l3}
 \ell_3(\eta)=d_\infty+\frac{\sum_{L\ne L_*}a_Ld_L-r}{a_\infty}\quad\text{in }\F_5 .
\end{equation}

Suppose that $m$ has no prime factor in $\cB$. The factorization reads
$m=r\prod_L\Nm(\eta_L)\cdot\Nm(\eta)^5$, the product including $L=L_*$, so no prime
factor of these norms lies in $\cB$. By the proof of Proposition~\ref{prop:g5-norms},
every prime ideal above such a prime has $d=0$, and so $\ell_3(\eta_L)=\ell_3(\eta)=0$
by additivity of $d$. Now~\eqref{eq:g5-boundary-l3} forces $r=0$ in $\F_5$,
contrary to $5\nmid r$. Hence $m$ has a prime factor in $\cB$. Since
Theorem~\ref{thm:A}(c) excludes the valuations from $1$ to $12$, this proves that
$S_0(G_5)$ is contained in the union in Theorem~\ref{thm:C}(b).

For the reverse inclusion, put
\[
 f_0=Y+Y^4-1-X,\qquad \Sigma'=\sum_{i,j=0}^4X^iY^j .
\]
The polynomial $f_0+\Sigma'$ is obtained by interchanging the variables in
$F_5$ of~\cite[\S3, (14)]{Panraksa2026Order25}; here we also need
$\Delta$ for its lifts.
Its augmentation is zero, and its values $\beta_L=f_0(\omega^{L_1},\omega^{L_2})$ are
\[
\begin{array}{c|cccccc}
 L&L_\infty&L_0&L_1&L_2&L_3&L_4\\ \hline
 \beta_L&\pi&\omega^4\pi^2&\omega^4\pi&
 \omega^4(1+\omega^2)\pi&-(1+\omega^3)\pi&-\pi
\end{array}
\]
The factors other than the powers of $\pi$ are global units, so
$\prod_L\Nm(\beta_L)=5^7$ and $M_1(f_0+t\Sigma')=25t\cdot5^7=5^9t$ for every
$t\in\Z$.

Consider the family
\begin{equation}\label{eq:g5-explicit-family}
 F_t=y+y^4-1-x+t\,(x+x^2+x^3+x^4+x^5)(1+y+y^2+y^3+y^4),\qquad t\in\Z,
\end{equation}
whose abelianization is $f_0+t\Sigma'$. To compute $\Delta(F_t)$ we use the matrices
of Lemma~\ref{lem:g5-singular} and put
\[
 A_0=\omega+\omega^4-1,\qquad B_0=\omega^2+\omega^3-1 ,
\]
so that $A_0+B_0=-3$ and $A_0B_0=1$. The matrix $\mathsf U+\mathsf U^4-I-\mathsf T$
has diagonal entries $1,A_0,B_0,B_0,A_0$ and determinant $1-\omega=\pi$. For the
auxiliary element
\[
 F_t^\circ=y+y^4-1-x+t\Bigl(\sum_{i=0}^4x^i\Bigr)\Bigl(\sum_{j=0}^4y^j\Bigr),
\]
all entries depending on $t$ lie in column zero, because
$\sum_{j=0}^4\mathsf U^j=5E_{00}$, where $E_{00}$ is the matrix unit. Expanding along
that column gives
\[
 \Delta(F_t^\circ)=\pi+5t\bigl[1+\omega(1+A_0+A_0B_0+A_0B_0^2)\bigr]=(1+5t)\pi ,
\]
because $1+A_0+A_0B_0+A_0B_0^2=2+A_0+B_0=-1$. Passing from $F_t^\circ$ to $F_t$ adds
$t(\mathsf T^5-I)\sum_j\mathsf U^j=-5t\pi E_{00}$. The $(0,0)$-cofactor does not
depend on $t$ and equals $A_0^2B_0^2=1$. Consequently $\Delta(F_t)=\pi$ for every
$t\in\Z$, and Lemma~\ref{lem:fact} gives
\[
 D_{G_5}(F_t)=5^9t\cdot5^5=5^{14}t .
\]
This gives $5^{14}\Z\subseteq S(G_5)$.

It remains to realize $5^{13}qr$ for $q\in\cB$ and $r$ prime to $5$; the case
$5\mid r$ is already covered. Choose a generator $\gamma\in\cO$ of a prime ideal of
norm $q$ and put $d=\ell_3(\gamma)\ne0$. By Lemma~\ref{lem:g5-units} we may multiply
$\gamma$ by a global unit so that its residue, $\ell_1$ and $\ell_2$ become $1$, $0$
and $0$. Its norm is still $q$, its invariant $\ell_3$ is still $d$, and
\[
 \gamma\equiv1+dw^3\equiv1-d\pi^3\pmod{\pi^4\cO}.
\]
Choose $a\in\{1,2,3,4\}$ representing $r/d$ in $\F_5$. We prescribe the augmentation
$c=5r$ and the six character values $\alpha_{L_a}=\beta_{L_a}\gamma$ and
$\alpha_L=\beta_L$ for $L\ne L_a$. Since $\beta_{L_a}/\pi\equiv a\pmod\pi$, the
change at $L_a$ satisfies
\[
 \xi_0:=\alpha_{L_a}-\beta_{L_a}=\beta_{L_a}(\gamma-1)\equiv-ad\pi^4\equiv5r
 \pmod{\pi^5\cO}.
\]
Apply Fourier inversion as in the proof of Lemma~\ref{lem:g5-reconstruct}, with
coefficients
\[
 b_{ij}=\frac1{25}\Bigl(c+\sum_L\Tr\bigl(\alpha_L\omega^{-iL_1-jL_2}\bigr)\Bigr),
 \qquad 0\le i,j<5 .
\]
Compared with the integral coefficients of $f_0$, each numerator changes by
$5r+\Tr(\xi_0\omega^{-ia-j})$. For every integer $h$ we have
$\xi_0\omega^h\equiv5r\pmod{\pi^5\cO}$, and hence
$\Tr(\xi_0\omega^h)\equiv20r\equiv-5r\pmod{25}$, because
$\Tr(\pi^5\cO)\subseteq25\Z$. Thus every $b_{ij}$ is an integer, and Fourier inversion
gives $f\in\Z[C_5^2]$ with exactly the prescribed values and augmentation. Their norm
product is $5^7q$, so $M_1(f)=5^8qr$.

The values $\alpha_L/\pi$ for $L\ne L_0$ have $\ell_3=0$ except on $L_a$, where the
residue is $a$ and $\ell_3$ is $d$, and $\alpha_{L_0}$ still has valuation $2$.
By~\eqref{eq:g5-boundary-trace}, $\theta(f)/5\equiv ad\equiv r\pmod5$. Since
$\alpha_\infty=\pi$, Lemma~\ref{lem:g5-divisible-scalar} gives
\[
 \Delta(F')\equiv\pi\pmod{\pi^5\cO}
\]
for every lift $F'\in\Z[G_5]$ of $f$. In particular, $\Delta(F')$ has valuation $1$,
so by Lemma~\ref{lem:g5-singular} the element $\pi$ lies in
$\nrd(\rho(F')+\pi\Lam_\p)$, and Proposition~\ref{prop:fibre} gives a lift $F$ of
$f$ with $\Delta(F)=\pi$. Then $D_{G_5}(F)=5^{13}qr$ by Lemma~\ref{lem:fact}. Here
$r$ may have either sign and any prime factors other than $5$, including $q$ itself,
so this realizes the whole union.

Finally, $11$ is the smallest prime congruent to $1$ modulo $5$, and $11\in\cB$ by
Lemma~\ref{lem:cubic-obstruction}. Hence a positive group determinant of valuation $13$
is at least $11\cdot5^{13}$, while $5^{14}=D_{G_5}(F_t)$ with $t=1$. This proves the remaining
assertions of Theorem~\ref{thm:C}(b). \qed

\section{Further questions}\label{sec:open}

Theorem~\ref{thm:A}(b) gives $T_3\subseteq S(\Gp)$ for every odd prime $p$. By
Lemma~\ref{lem:L1} and Proposition~\ref{prop:fibre}, once $\Fbar$ is fixed with
$p\nmid F(1)$, the factor $\Delta$ can be any element of $\cO\cap\Delta(F)(1+\p^p\cO_\p)$; what
remains is to decide which combinations of linear character values occur.

\begin{problem}\label{prob:suff}
Determine $S^*(\Gp)$ for $p\ge7$. Which conditions must be added to
Proposition~\ref{prop:nec} to obtain a sufficient criterion?
\end{problem}

At $p=5$ the global units adjust the residue, $\ell_1$ and $\ell_2$ but not $\ell_3$
or the norm modulo $25$, and this produces the extra condition in
Theorem~\ref{thm:C}. For $p\ge7$ more coefficients of the logarithm lie outside the
reach of the global units~\cite[\S\S3--5]{Panraksa2026Bernoulli}, which leads to
further conditions. The set
$S_0(\Gp)$ is also unknown for $p\ge7$; in particular, we do not know whether
$p^{2p+3}$ itself is a group determinant of $\Gp$.

\section*{Statements and Declarations}

\textbf{Funding.} This work was supported by Mahidol University International
College (MUIC).

\smallskip
\textbf{Competing interests.} The author has no competing interests to declare.

\smallskip
\textbf{Data availability.} No data sets were generated or analysed. The Magma and
Python scripts used to check the explicit examples are available from the author on
reasonable request.

\bibliographystyle{amsplain}
\bibliography{references}

\end{document}